# Approximation of the relaxed Euler system at the level of the compressible Navier-Stokes equations

Hao Chen*[a], Tao Wang†[a], Junhao Zhang‡[b], and Huijiang Zhao§[a]

[a]Wuhan University, School of Mathematics and Statistics, Wuhan 430072, China
[b]University of Warsaw, Institute of Applied Mathematics and Mechanics, Banacha 2, 02-097 Warszawa, Poland

**Abstract**

We study the global well-posedness of the Cauchy problem in $\mathbb{R}^d$, $d \geq 2$, for the relaxed compressible Euler system with a Maxwell-type constitutive law for the stress tensor. We show that the corresponding compressible Navier-Stokes solutions serve as both relaxation limits and large-time asymptotic profiles of the global relaxed Euler solutions constructed here. Given a Navier-Stokes background satisfying suitable boundedness and decay assumptions, we construct global smooth relaxed Euler solutions for each fixed $\alpha \in (0,1]$ and sufficiently small $\varepsilon > 0$. Here, the background need not be small. After subtracting the stress initial layer, the approximation error is $O(\varepsilon^\alpha)$ uniformly in time. Under suitable additional conditions on the initial perturbations, this error decays in $L^2$ faster than the decay of standard Navier-Stokes background, uniformly in $\varepsilon$, identifying the background solution as a large-time asymptotic profile. As an application, we obtain the corresponding small-data results using classical Navier-Stokes theory to verify the background assumptions.

**Keywords:** Relaxed compressible Euler system; compressible Navier-Stokes equations; global well-posedness; uniform-in-time approximation error; asymptotic profiles.

# 1 Introduction and main results

To allow the stress tensor of a viscous compressible fluid to relax toward equilibrium at a finite rate instead of responding instantaneously, Maxwell [22] proposed a relaxation-type constitutive law governing its evolution. Accordingly, the dynamics of an isentropic compressible fluid are described by the following relaxed Euler system:

$$\partial_t \rho^\varepsilon + \mathrm{div}(\rho^\varepsilon u^\varepsilon) = 0,$$

*Email: chenhao215@outlook.com
†Email: tao.wang@whu.edu.cn
‡Email: jhzhang@mimuw.edu.pl
§Email: hhjjzhao@whu.edu.cn

$$\begin{aligned} \partial_t(\rho^\varepsilon u^\varepsilon) + \operatorname{div}(\rho^\varepsilon u^\varepsilon \otimes u^\varepsilon) + \nabla p(\rho^\varepsilon) &= -\operatorname{div}\pi^\varepsilon, \\ \varepsilon\partial_t\pi^\varepsilon + \pi^\varepsilon &= -\nu\sigma(u^\varepsilon) - \lambda(\operatorname{div} u^\varepsilon)I_d, \end{aligned} \tag{1.1}$$

where $x \in \mathbb{R}^d$ is the spatial variable, $t \in \mathbb{R}^+$ is the time variable, and $d \geq 2$ is the spatial dimension. Here, $\rho^\varepsilon(t,x) \geq 0$, $u^\varepsilon(t,x) \in \mathbb{R}^d$ and $\pi^\varepsilon(t,x) \in \mathbb{R}^{d\times d}$ denote the fluid density, the velocity field and the symmetric viscous stress tensor. The pressure $p = p(\rho^\varepsilon)$ is a smooth function of the density satisfying $p'(\rho^\varepsilon) > 0$, for all $\rho^\varepsilon > 0$. We define the trace-free symmetric tensor by

$$\sigma(u^\varepsilon) = \nabla u^\varepsilon + (\nabla u^\varepsilon)^T - \frac{2}{d}(\operatorname{div} u^\varepsilon)I_d,$$

where $I_d$ denotes the $d \times d$ identity matrix. The constant $\varepsilon > 0$ is the relaxation time, while $\nu$ and $\lambda$ are the viscosity coefficients satisfying

$$\nu > 0, \qquad \lambda > 0, \qquad d \geq 2.$$

We first address the global well-posedness of the Cauchy problem for (1.1) with the initial data

$$\begin{aligned} (\rho^\varepsilon, u^\varepsilon, \pi^\varepsilon)(0,x) &= (\rho_0^\varepsilon, u_0^\varepsilon, \pi_0^\varepsilon)(x), \quad x \in \mathbb{R}^d, \\ \pi_0^\varepsilon(x) &= \left(\pi_0^\varepsilon(x)\right)^T, \quad x \in \mathbb{R}^d, \\ \lim_{|x|\to\infty}(\rho_0^\varepsilon, u_0^\varepsilon, \pi_0^\varepsilon)(x) &= (1,0,0). \end{aligned} \tag{1.2}$$

We recall the well-posedness theory for (1.1) and related stress-relaxation models. Since (1.1) admits a strictly convex entropy [24], the theory of symmetrizable hyperbolic systems [1, 4, 5, 15, 16, 19] yields a unique local smooth solution with $(\rho^\varepsilon(t,x) - 1, u^\varepsilon(t,x), \pi^\varepsilon(t,x)) \in C\left([0,T); H^s(\mathbb{R}^d)\right) \cap C^1\left([0,T); H^{s-1}(\mathbb{R}^d)\right)$. For a related three-dimensional isentropic model, Yong [33] introduced separate relaxation equations for the deviatoric and volumetric stresses and established an entropy-dissipation structure that makes the system symmetrizable hyperbolic.

However, local smooth solutions cannot always be extended for all time. Sideris [28] established finite-time singularity formation for nonlinear hyperbolic conservation laws under suitable structural and initial-data assumptions. For compressible flow with revised Maxwell's law, Wang and Hu [32] proved finite-time breakdown of $C^1$ solutions for a suitable class of large initial data, using finite propagation and a Sideris-type argument. Nevertheless, for $s > d/2+1$, Peng [24] established global existence of smooth solutions to (1.1) for small perturbations of the constant equilibrium $(1,0,0)$.

Related results are available for non-isentropic models. In two and three dimensions, Hu and Racke [9] proved local well-posedness and global existence for small Sobolev perturbations for a system combining revised Maxwell stress relaxation with Fourier heat conduction. They subsequently considered simultaneous relaxation of the viscous stress and heat flux [10], treating two special cases: small data global existence when the shear viscosity remains positive, and finite-time blow-up for a class of large data in the corresponding purely hyperbolic case. For three-dimensional compressible Navier-Stokes equations with the Cattaneo heat conduction, global existence and optimal time-decay estimates for small initial perturbations were established in [17, 30].

We further investigate the approximation of relaxed Euler solutions by Navier-Stokes flows, both in the relaxation limit $\varepsilon \to 0$ and in the large-time regime $t \to \infty$. The analysis relies on the existence and large-time behavior of solutions to the Navier-Stokes equations:

$$
\begin{aligned}
\partial_t \rho + \operatorname{div}(\rho u) &= 0, \\
\partial_t(\rho u) + \operatorname{div}(\rho u \otimes u) + \nabla p(\rho) &= \nu \Delta u + \left(\nu + \lambda - \frac{2\nu}{d}\right) \nabla(\operatorname{div} u).
\end{aligned} \tag{1.3}
$$

Here, the associated viscous stress tensor is given by

$$
\pi = -\nu\sigma(u) - \lambda(\operatorname{div} u) I_d.
$$

We supplement the system (1.3) with the initial data

$$
\begin{aligned}
(\rho, u)(0, x) &= (\rho_0, u_0)(x), \\
\lim_{|x|\to\infty} (\rho_0, u_0)(x) &= (1, 0).
\end{aligned} \tag{1.4}
$$

Classical contributions to the local theory of smooth compressible flow include those of Serrin [27] and Nash [23]. Matsumura and Nishida [20, 21] developed the global theory for small smooth perturbations of constant equilibrium states. In a different regularity framework, Danchin [3] established global well-posedness near equilibrium in critical Besov spaces in dimensions $d \geq 2$.

The large-time behavior of small solutions has been studied through both Green function analysis and energy methods. Hoff and Zumbrun [7, 8] described multidimensional diffusion waves and derived pointwise decay estimates. Liu and Wang [18] obtained pointwise estimates for the asymptotic diffusion waves in odd space dimensions. For the three-dimensional isentropic system, Guo and Wang [6] developed an energy method based on homogeneous negative Sobolev norms to recover time decay, including decay of higher spatial derivatives, without using linear decay estimates. These works give more detailed information than decay in Sobolev norms alone.

For large initial data, a separate line of work concerns two-dimensional systems with density-dependent bulk viscosity. This includes the global existence theory of Vaigant and Kazhikhov [31], the Cauchy problem with non-vacuum far-field states studied by Jiu, Wang and Xin [12, 13], and the vacuum far-field problem treated by Huang and Li [11]. These results allow large data under suitable assumptions on the viscosity law and initial state. Their viscosity laws or far-field conditions differ from those of the constant-viscosity, non-vacuum problem considered here.

In view of the existence and decay results recalled above, we assume that $(\rho_0(x)-1, u_0(x)) \in H^{s_1}(\mathbb{R}^d) \cap L^1(\mathbb{R}^d)$ and that the Cauchy problem (1.3), (1.4) admits a global solution $(\rho(t,x), u(t,x))$ satisfying the following regularity, boundedness and decay conditions:

- For $s_1 > d + [d/2] + 11$, assume that $(\rho(t,x) - 1, u(t,x)) \in C\left([0,\infty); H^{s_1}(\mathbb{R}^d)\right) \cap C^1\left([0,\infty); H^{s_1-2}(\mathbb{R}^d)\right)$, and there exist positive constants $\underline{\rho}$, $\bar{\rho}$ and $\bar{u}$, independent of $t$ and $x$, such that
$$
\underline{\rho} \leq \rho(t,x) \leq \bar{\rho}, \qquad |u(t,x)| \leq \bar{u}, \qquad (t,x) \in [0,\infty) \times \mathbb{R}^d. \tag{H$_1$}
$$
- For $d \geq 2$, $0 \leq k \leq \frac{s_1 - [d/2] - 3}{2}$, $\|(\rho_0 - 1, u_0)\|_{H^{s_1} \cap L^1(\mathbb{R}^d)} \lesssim 1$ and $(\rho(t,x), u(t,x))$ satisfies the decay estimates
$$
\begin{aligned}
&\left\| D^k(\rho(t) - 1, u(t)) \right\|_{L^p(\mathbb{R}^d)} \\
\lesssim &\|(\rho_0 - 1, u_0)\|_{H^{s_1} \cap L^1(\mathbb{R}^d)} (1+t)^{-\frac{d}{2}\left(1 - \frac{1}{p}\right) - \frac{k}{2}}, \quad 2 \leq p \leq \infty, \quad t \geq 0.
\end{aligned} \tag{H$_2$}
$$

Here, $[a]$ denotes the greatest integer not exceeding $a$. It is clear that, under suitable smallness assumptions on the initial data in $H^{s_1}(\mathbb{R}^d)\cap L^1(\mathbb{R}^d)$, global smooth solutions of the compressible Navier-Stokes equations satisfy $(\mathrm{H}_1)$ and $(\mathrm{H}_2)$ (see [7, 8, 14, 18, 20, 21]).

We next review the convergence of relaxed Euler solutions to Navier-Stokes flows. Yong [33] justified the Newtonian limit for his three-dimensional isentropic Maxwell model on the lifespan of a smooth Navier-Stokes solution. Hu and Racke [9] obtained finite-time relaxation estimates for the related non-isentropic model with Fourier heat conduction. Peng [24] studied both tensor and vector formulations of relaxed Euler systems, establishing quantitative convergence on fixed time intervals with initial-layer corrections. For initial data near constant equilibrium states, he also proved uniform global existence and global-in-time convergence by compactness.

Global relaxation limits have also been established for other constitutive laws. Peng and Zhao [25] considered a one-dimensional full compressible model with Oldroyd-type stress relaxation and Cattaneo-Christov heat conduction, proving global existence and convergence near equilibrium. For the corresponding periodic problem, they obtained global convergence rates by combining stream functions with energy estimates [26]. For the Navier-Stokes-Cattaneo-Christov system in dimensions $d \geq 3$, Crin-Barat, Kawashima and Xu [2] established global well-posedness near equilibrium in a critical regularity framework, sharp time decay and global-in-time quantitative relaxation estimates toward the Navier-Stokes-Fourier system. This approximation retains Newtonian viscosity and relaxes the heat flux.

The present work concerns stress relaxation around a prescribed Navier-Stokes background satisfying $(\mathrm{H}_1)$ and $(\mathrm{H}_2)$, without imposing any smallness condition on its initial data.

To state the main results, write $\pi_0(x) = \pi(0,x)$, let $\alpha \geq 0$, set $m^\varepsilon(t,x) = \rho^\varepsilon(t,x)u^\varepsilon(t,x)$ and $m(t,x) = \rho(t,x)u(t,x)$ and introduce

$$\begin{aligned} \rho^\varepsilon(t,x) &= \rho(t,x) + \varepsilon^\alpha z^\varepsilon(t,x), \\ u^\varepsilon(t,x) &= u(t,x) + \varepsilon^\alpha w^\varepsilon(t,x), \\ \pi^\varepsilon(t,x) &= \pi(t,x) + \varepsilon^\alpha \chi^\varepsilon(t,x), \\ m^\varepsilon(t,x) &= m(t,x) + \varepsilon^\alpha n^\varepsilon(t,x). \end{aligned} \tag{1.5}$$

A direct computation gives

$$n^\varepsilon = \rho w^\varepsilon + u z^\varepsilon + \varepsilon^\alpha z^\varepsilon w^\varepsilon. \tag{1.6}$$

The Cauchy problem (1.1), (1.2) can therefore be rewritten as

$$\begin{cases} z_t^\varepsilon + u^\varepsilon \cdot \nabla z^\varepsilon + \rho^\varepsilon \operatorname{div} w^\varepsilon = \widetilde{F}_1, \\ w_t^\varepsilon + u^\varepsilon \cdot \nabla w^\varepsilon + h'(\rho^\varepsilon)\nabla z^\varepsilon + \dfrac{\operatorname{div}\chi^\varepsilon}{\rho^\varepsilon} = \widetilde{F}_2, \\ \varepsilon \chi_t^\varepsilon + \chi^\varepsilon + \nu\sigma(w^\varepsilon) + \lambda(\operatorname{div} w^\varepsilon) I_d + \varepsilon^{1-\alpha}\pi_t = 0, \end{cases} \tag{1.7}$$

with initial data

$$\begin{aligned} (z^\varepsilon, w^\varepsilon, \chi^\varepsilon)|_{t=0} &= (z_0^\varepsilon, w_0^\varepsilon, \chi_0^\varepsilon) \\ &:= \varepsilon^{-\alpha}(\rho_0^\varepsilon - \rho_0, u_0^\varepsilon - u_0, \pi_0^\varepsilon - \pi_0) \\ &\longrightarrow (0,0,0), \qquad |x| \to \infty, \end{aligned} \tag{1.8}$$

where

$$\widetilde{F}_1 = -z^\varepsilon \operatorname{div} u - w^\varepsilon \cdot \nabla\rho, \qquad h'(\rho^\varepsilon) = \frac{p'(\rho^\varepsilon)}{\rho^\varepsilon},$$

$$\widetilde{F}_2 = -\frac{z^\varepsilon \nabla\rho \int_0^1 p''(\rho+\theta\varepsilon^\alpha z^\varepsilon)\,d\theta}{\rho^\varepsilon} + \frac{p'(\rho) z^\varepsilon \nabla\rho}{\rho\rho^\varepsilon} + \frac{z^\varepsilon \operatorname{div}\pi}{\rho\rho^\varepsilon} - w^\varepsilon\cdot\nabla u.$$

Without loss of generality, we assume that $h(1) = 0$. Set $n_0^\varepsilon(x) = n^\varepsilon(0,x)$, $m_0^\varepsilon(x) = m^\varepsilon(0,x)$ and $m_0(x) = m(0,x)$. Our first theorem establishes global well-posedness and uniform-in-time relaxation convergence.

**Theorem 1.1.** *Assume that*

- $(\rho_0(x)-1, u_0(x)) \in H^{s_1}(\mathbb{R}^d)\cap L^1(\mathbb{R}^d)$, *and the Cauchy problem* (1.3), (1.4) *for the compressible Navier-Stokes equations admits a unique global solution satisfying* $(\mathrm{H}_1)$ *and* $(\mathrm{H}_2)$,
- $d\geq 2$, $s\in\mathbb{N}$ *with* $\frac{d}{2}+1 < s \leq \frac{s_1-9-\left[\frac{d}{2}\right]}{2}$, *and* $\alpha\in(0,1]$,
- $(z_0^\varepsilon(x), w_0^\varepsilon(x), \chi_0^\varepsilon(x)) \in H^s(\mathbb{R}^d)$ *and*

$$\delta_0 := \|(z_0^\varepsilon, w_0^\varepsilon, \sqrt{\varepsilon}\,\chi_0^\varepsilon)\|_{H^s(\mathbb{R}^d)} \lesssim 1. \tag{1.9}$$

*Then there exists a sufficiently small* $\varepsilon_0 > 0$, *depending on the fixed* $\alpha$, *a uniform upper bound for* $\delta_0$ *in* (1.9), *for all* $0<\varepsilon\leq\varepsilon_0$, *the Cauchy problem* (1.7), (1.8) *admits a unique global smooth solution* $(z^\varepsilon(t,x), w^\varepsilon(t,x), \chi^\varepsilon(t,x)) \in C\left([0,\infty); H^s(\mathbb{R}^d)\right)\cap C^1\left([0,\infty); H^{s-1}(\mathbb{R}^d)\right)$ *satisfying*

$$\begin{aligned}&\|(z^\varepsilon(t), w^\varepsilon(t), \sqrt{\varepsilon}\chi^\varepsilon(t))\|_{H^s(\mathbb{R}^d)}^2 + \int_0^t \left(\|\nabla(z^\varepsilon(\tau), w^\varepsilon(\tau))\|_{H^{s-1}(\mathbb{R}^d)}^2 + \|\chi^\varepsilon(\tau)\|_{H^s(\mathbb{R}^d)}^2\right)d\tau\\ \lesssim&\|(z_0^\varepsilon, w_0^\varepsilon, \sqrt{\varepsilon}\chi_0^\varepsilon)\|_{H^s(\mathbb{R}^d)}^2 + \varepsilon^{2-2\alpha},\end{aligned} \tag{1.10}$$

*and*

$$\|\chi^\varepsilon(t) - e^{-\frac{t}{\varepsilon}}\chi_0^\varepsilon\|_{H^{s-1}(\mathbb{R}^d)} \lesssim \|(z_0^\varepsilon, w_0^\varepsilon, \sqrt{\varepsilon}\,\chi_0^\varepsilon)\|_{H^s(\mathbb{R}^d)} + \varepsilon^{1-\alpha}. \tag{1.11}$$

*Consequently, by* (1.5), (1.10) *and* (1.11), *the Cauchy problem* (1.1), (1.2) *for the relaxed Euler system admits a unique global smooth solution* $(\rho^\varepsilon(t,x), u^\varepsilon(t,x), \pi^\varepsilon(t,x))$. *Moreover, the density and velocity errors, together with the stress error after subtraction of its initial layer, satisfy the uniform-in-time estimate*

$$\begin{aligned}&\|(\rho^\varepsilon(t)-\rho(t), u^\varepsilon(t)-u(t))\|_{H^s(\mathbb{R}^d)}\\ &\quad + \|\pi^\varepsilon(t) - \pi(t) - e^{-\frac{t}{\varepsilon}}(\pi_0^\varepsilon - \pi_0)\|_{H^{s-1}(\mathbb{R}^d)} \lesssim \varepsilon^\alpha, \qquad t\geq 0.\end{aligned} \tag{1.12}$$

*Remark* 1.1. In the proof of Theorem 1.1, we only use the uniform boundedness estimates $(\mathrm{H}_1)$ and the time decay estimates $(\mathrm{H}_2)$ of the global solution $(\rho(t,x), u(t,x))$ to the Cauchy problem (1.3), (1.4) for the compressible Navier-Stokes equations. No smallness condition is imposed on this solution. Therefore, if a global solution with large initial data satisfying $(\mathrm{H}_1)$ and $(\mathrm{H}_2)$ is available, then for sufficiently small relaxation time $\varepsilon > 0$, one can construct a global smooth solution to the Cauchy problem (1.1), (1.2) for the relaxed Euler system near this background solution.

*Remark* 1.2. In Theorem 1.1, for $0<\alpha\leq 1$, we only require the relaxation time $\varepsilon>0$ to be sufficiently small. No smallness assumption is needed for the initial data $(z_0^\varepsilon(x), w_0^\varepsilon(x), \chi_0^\varepsilon(x))$.

*Remark* 1.3. By (1.9), (1.10) and (1.12), the uncorrected stress error satisfies

$$\|\pi^\varepsilon(t)-\pi(t)\|_{H^{s-1}(\mathbb{R}^d)}\lesssim \varepsilon^\alpha+\varepsilon^{\alpha-\frac12}e^{-\frac{t}{\varepsilon}},\qquad t\geq 0.$$

Consequently, for every fixed $t_0>0$, the stress converges at rate $O(\varepsilon^\alpha)$ uniformly for $t\geq t_0$, with a constant that may depend on $t_0$. For $\alpha>\frac12$, convergence is also uniform for $t\geq 0$, at rate $O(\varepsilon^{\alpha-\frac12})$. When $0<\alpha\leq\frac12$, the initial-error condition (1.9) alone does not ensure convergence of the initial stress.

The following theorem identifies solutions of the compressible Navier-Stokes system as large-time asymptotic profiles of the corresponding solutions of the relaxed compressible Euler system.

**Theorem 1.2.** *In addition to the assumptions of Theorem 1.1, assume that*

- $s\geq\left[\frac d2\right]+4$,
- *the initial data* $(z_0^\varepsilon(x),n_0^\varepsilon(x),\chi_0^\varepsilon(x))$ *satisfy*

$$\int_{\mathbb{R}^d}(z_0^\varepsilon(x),n_0^\varepsilon(x))dx=(0,0),\quad \left\|(1+|x|)\big(z_0^\varepsilon(x),n_0^\varepsilon(x)\big)\right\|_{L^1(\mathbb{R}^d)}+\sqrt{\varepsilon}\|\chi_0^\varepsilon(x)\|_{L^1(\mathbb{R}^d)}\lesssim 1. \tag{1.13}$$

*Then there exists a sufficiently small positive constant* $\varepsilon_1\leq\varepsilon_0$ *such that for all* $0<\varepsilon\leq\varepsilon_1$ *and* $t\geq 0$, *the global solution* $(z^\varepsilon(t,x),w^\varepsilon(t,x),\chi^\varepsilon(t,x))$ *constructed in Theorem 1.1 satisfies the following time decay estimates uniformly with respect to the relaxation time* $\varepsilon$:

$$\|D^k(z^\varepsilon(t),w^\varepsilon(t))\|_{L^2(\mathbb{R}^d)}\lesssim\begin{cases}(1+t)^{-\frac{d+2k+2}{4}}, & k=0,1,\dots,\left[\frac d2\right]+1,\\ (1+t)^{-\frac{d+2[\frac d2]+4}{4}}, & k=\left[\frac d2\right]+2,\end{cases} \tag{1.14}$$

*and*

$$\|D^k(\chi^\varepsilon(t)-e^{-\frac{t}{\varepsilon}}\chi_0^\varepsilon)\|_{L^2(\mathbb{R}^d)}\lesssim\begin{cases}(1+t)^{-\frac{d+2k+4}{4}}, & k=0,1,\dots,\left[\frac d2\right],\\ (1+t)^{-\frac{d+2[\frac d2]+4}{4}}, & k=\left[\frac d2\right]+1.\end{cases} \tag{1.15}$$

*Remark* 1.4. The decay rates in Theorem 1.2 for derivatives below the highest orders are optimal for the linearized system for each $\varepsilon>0$. They are half an order faster than the standard decay rates for Navier-Stokes solutions, owing to low-frequency cancellation. Indeed, the zero-mean conditions in (1.13) eliminate the constant term in the Fourier expansion, while the finite first-moment bound gives $(\widehat{z_0^\varepsilon}(\xi),\widehat{n_0^\varepsilon}(\xi))=O(|\xi|)$ near the origin. The quadratic low-frequency dissipation then converts this additional factor of $|\xi|$ into a gain of $(1+t)^{-1/2}$ in the lower-order decay estimates.

*Remark* 1.5. The half-order loss at the highest derivatives arises because the present $H^1$ estimate controls derivatives of orders $[d/2]+1$ and $[d/2]+2$ with the same decay rate. The highest-order corrected stress estimate inherits this loss through (2.4), which involves one spatial derivative of the velocity.

*Remark* 1.6. If $(\rho_0^\varepsilon(x),u_0^\varepsilon(x),\pi_0^\varepsilon(x))=(\rho_0(x),u_0(x),\pi_0(x))$, then $(z_0^\varepsilon(x),n_0^\varepsilon(x),\chi_0^\varepsilon(x))=0$. Hence the initial-error conditions in (1.13) hold automatically and the stress initial layer vanishes. Under the remaining assumptions of Theorem 1.2, Theorems 1.1 and 1.2 yield a unique global smooth solution for sufficiently small $\varepsilon > 0$. Combining (1.14) and (1.15) with ($\mathrm{H}_2$) gives the standard decay estimates, uniformly in $\varepsilon$,

$$\|D^k(\rho^\varepsilon(t)-1, u^\varepsilon(t))\|_{L^2(\mathbb{R}^d)} \lesssim (1+t)^{-\frac{d+2k}{4}}, \quad 0 \le k \le [d/2]+1,$$
$$\|D^k\pi^\varepsilon(t)\|_{L^2(\mathbb{R}^d)} \lesssim (1+t)^{-\frac{d+2k+2}{4}}, \quad 0 \le k \le [d/2].$$

Since the small-data theory [3, 6, 7, 18, 20, 21] provides a global Navier-Stokes background satisfying ($\mathrm{H}_1$) and ($\mathrm{H}_2$), we obtain the following global existence and decay theorem by adapting the arguments used to prove Theorems 1.1 and 1.2. The result holds uniformly for $\varepsilon \in (0,1]$.

**Theorem 1.3.** *Assume that*

- $d \ge 2$ *and* $s, s_1 \in \mathbb{N}$ *satisfy* $\frac{d}{2}+1 < s \le \frac{s_1-[\frac{d}{2}]-9}{2}$,
- *the initial data* $(\rho_0(x)-1, u_0(x)) \in H^{s_1}(\mathbb{R}^d)\cap L^1(\mathbb{R}^d)$ *and* $(\rho_0^\varepsilon(x)-\rho_0(x), u_0^\varepsilon(x)-u_0(x), \pi_0^\varepsilon(x)-\pi_0(x)) \in H^s(\mathbb{R}^d)$ *satisfy*

$$\big\|\big(\rho_0^\varepsilon-\rho_0,\, u_0^\varepsilon-u_0,\, \sqrt{\varepsilon}(\pi_0^\varepsilon-\pi_0)\big)\big\|_{H^s(\mathbb{R}^d)} \le \varepsilon^\alpha \delta_1,$$
$$\|(\rho_0-1, u_0)\|_{H^{s_1}\cap L^1(\mathbb{R}^d)} \le \delta_2.$$

*There exist constants* $\tilde{\epsilon}_0 > 0$, *independent of* $\varepsilon \in (0,1]$ *and* $\alpha \ge 0$, *such that, whenever* $\delta_1+\delta_2 \le \tilde{\epsilon}_0$, *the Cauchy problem for the compressible Navier-Stokes equations* (1.3), (1.4) *admits a unique global smooth solution* $(\rho(t,x), u(t,x), \pi(t,x))$ *satisfying* ($\mathrm{H}_1$) *and* ($\mathrm{H}_2$). *The Cauchy problem for the relaxed Euler system* (1.1), (1.2) *admits a unique global smooth solution* $(\rho^\varepsilon(t,x), u^\varepsilon(t,x), \pi^\varepsilon(t,x))$ *satisfying* $(\rho^\varepsilon(t,x)-1, u^\varepsilon(t,x), \pi^\varepsilon(t,x)) \in C\left([0,\infty); H^s(\mathbb{R}^d)\right) \cap C^1\left([0,\infty); H^{s-1}(\mathbb{R}^d)\right)$. *Furthermore, this solution satisfies the following energy estimate:*

$$\begin{aligned}&\big\|\big(\rho^\varepsilon(t)-\rho(t),\, u^\varepsilon(t)-u(t),\, \sqrt{\varepsilon}(\pi^\varepsilon(t)-\pi(t))\big)\big\|_{H^s(\mathbb{R}^d)}^2\\ &+\int_0^t \Big(\|\nabla(\rho^\varepsilon(\tau)-\rho(\tau), u^\varepsilon(\tau)-u(\tau))\|_{H^{s-1}(\mathbb{R}^d)}^2 + \|(\pi^\varepsilon(\tau)-\pi(\tau))\|_{H^s(\mathbb{R}^d)}^2\Big)\, d\tau\\ &\lesssim \varepsilon^{2\alpha}\delta_1^2+\varepsilon^2\delta_2^2,\end{aligned}$$

*and*

$$\begin{aligned}&\|(\rho^\varepsilon(t)-\rho(t), u^\varepsilon(t)-u(t))\|_{H^s(\mathbb{R}^d)}\\ &\quad+\Big\|\pi^\varepsilon(t)-\pi(t)-e^{-\frac{t}{\varepsilon}}(\pi_0^\varepsilon-\pi_0)\Big\|_{H^{s-1}(\mathbb{R}^d)} \lesssim \varepsilon^\alpha\delta_1+\varepsilon\delta_2.\end{aligned}$$

*For the time decay estimates, assume further that*

- $s \ge \left[\frac{d}{2}\right]+4$,
- *the initial data satisfy*

$$\int_{\mathbb{R}^d}\big(\rho_0^\varepsilon-\rho_0,\, m_0^\varepsilon-m_0\big)\, dx = (0,0),$$
$$\big\|(1+|x|)\big(\rho_0^\varepsilon-\rho_0,\, m_0^\varepsilon-m_0\big)\big\|_{L^1(\mathbb{R}^d)} + \sqrt{\varepsilon}\|\pi_0^\varepsilon-\pi_0\|_{L^1(\mathbb{R}^d)} \lesssim \varepsilon^\alpha.$$

*There exists a constant $0 < \tilde{\epsilon}_1 \le \tilde{\epsilon}_0$, independent of $\varepsilon$ and $\alpha \ge 0$, such that, whenever $\delta_1 + \delta_2 \le \tilde{\epsilon}_1$, the global solution satisfies the following time decay estimates uniformly with respect to the relaxation time $\varepsilon$:*

$$\|D^k(\rho^\varepsilon(t) - \rho(t), u^\varepsilon(t) - u(t))\|_{L^2(\mathbb{R}^d)} \lesssim \begin{cases} \varepsilon^{\min\{1,\alpha\}}(1+t)^{-\frac{d+2k+2}{4}}, & 0 \le k \le [d/2]+1, \\ \varepsilon^{\min\{1,\alpha\}}(1+t)^{-\frac{d+2[d/2]+4}{4}}, & k = [d/2]+2, \end{cases}$$

*and*

$$\left\|D^k\left(\pi^\varepsilon(t) - \pi(t) - e^{-\frac{t}{\varepsilon}}(\pi_0^\varepsilon - \pi_0)\right)\right\|_{L^2(\mathbb{R}^d)} \lesssim \begin{cases} \varepsilon^{\min\{1,\alpha\}}(1+t)^{-\frac{d+2k+4}{4}}, & 0 \le k \le [d/2], \\ \varepsilon^{\min\{1,\alpha\}}(1+t)^{-\frac{d+2[d/2]+4}{4}}, & k = [d/2]+1. \end{cases}$$

We outline the main difficulties and ideas in the proofs of Theorems 1.1 and 1.2.

For Theorem 1.1, the main difficulties are to allow a Navier-Stokes background that need not be small and to recover dissipation for the density and velocity errors. We construct the energy around the prescribed background using the relative entropy

$$\begin{aligned}\eta(t,x) = {} & \frac{\rho^\varepsilon|u^\varepsilon|^2}{2} - \frac{\rho|u|^2}{2} + H(\rho^\varepsilon) - H(\rho) - \left(h(\rho) - \frac{|u|^2}{2}\right)(\rho^\varepsilon - \rho) - u\cdot(\rho^\varepsilon u^\varepsilon - \rho u) \\ & + \frac{\varepsilon}{4\nu}|\pi^\varepsilon - \pi|^2 + \left(\frac{\varepsilon}{2d^2\lambda} - \frac{\varepsilon}{4\nu d}\right)|\mathrm{tr}(\pi^\varepsilon - \pi)|^2,\end{aligned}$$

where $H'(\rho) = h(\rho)$ and hence $H''(\rho) = p'(\rho)/\rho > 0$. The identity (2.8) gives the basic energy estimate for $(z^\varepsilon(t,x), w^\varepsilon(t,x), \chi^\varepsilon(t,x))$. For higher derivatives, we use the symmetrizer $B_0(\rho^\varepsilon) = \mathrm{diag}(h'(\rho^\varepsilon), \rho^\varepsilon I_d)$ and combine it with the quadratic stress energy. Since only the stress is directly dissipative, we use cross estimates from the momentum and Maxwell equations to recover the dissipation of $\nabla z^\varepsilon$ and $\nabla w^\varepsilon$.

To close these energy and dissipation estimates, we bound the background coupling terms by a small multiple of the dissipation and energy terms with time-integrable coefficients, using ($H_2$). The purely nonlinear perturbation terms carry the factor $\varepsilon^\alpha$ and can be absorbed by choosing $\varepsilon$ sufficiently small for fixed $\alpha$. Gronwall's inequality then yields the uniform energy estimate (1.10). Using this bound in the Maxwell representation (2.4), we obtain the stress estimate after subtraction of the initial layer, leading to the uniform-in-time error bound (1.12).

For Theorem 1.2, the main difficulty is to retain the low-frequency cancellation while controlling the background forcing uniformly in $\varepsilon$. We work with the density and momentum errors $(z^\varepsilon, n^\varepsilon)$ to preserve the conservative structure, and use Fourier energy estimates together with frequency and time decompositions. The coupling in (3.5) yields a frequency-dependent energy satisfying

$$G(t,\xi)^2 \sim |\xi|^{2k-2}(1+|\xi|^2)\left|\left(\widehat{z^\varepsilon}(t,\xi), \widehat{n^\varepsilon}(t,\xi)\right)\right|^2 + \varepsilon|\xi|^{2k}|\widehat{\chi^\varepsilon}(t,\xi)|^2.$$

To derive the desired decay estimates using the Fourier energy $G(t,\xi)^2$, we need low-frequency control of the background source $\varepsilon^{1-\alpha}\pi_t$. The Navier-Stokes equations express the linear part of $\pi_t$ in terms of spatial derivatives of $(\rho - 1, u)$. However, our assumptions do not directly provide a uniform-in-time $L^1$ bound for this background perturbation. This leads us to seek low-frequency control without propagating the spatial $L^1$ norm. Homogeneous negative Sobolev norms provide one way to quantify such information, as illustrated by the decay method of Guo and Wang [6]. For the present source estimate, uniform bounds in $\dot{H}^{-r}(\mathbb{R}^d)$ with $0 \le r < d/2$ alone do not yield the endpoint decay associated with $L^1$ data. We obtain the required endpoint control by estimating the Fourier transform directly in $L^\infty_\xi$. More precisely, exploiting the structure of (1.3), we establish the uniform bound in Lemma 3.2:

$$\sup_{t\ge 0}\big\|(\widehat{\rho(t)-1}, \widehat{u}(t))\big\|_{L^\infty_\xi} \lesssim 1.$$

Writing the background equations as a constant-coefficient linear system, we show that the nonlinear source terms belong to $L^1([0,\infty); L^1(\mathbb{R}^d))$. The initial $L^1$ bound and a Fourier energy estimate then yield the bound above. This bound is used to control the early-time low-frequency linear contribution of $\pi_t$ in (3.27), while $(\mathrm{H}_2)$ handles the other contributions to that estimate.

In the Fourier energy identity, we must also estimate the terms involving $F_1$ and $F_2$ in (3.4) while preserving the low-frequency weights and maintaining uniformity in $\varepsilon$. This requires different treatments of the two source terms. Direct absorption of the $F_1$ term into the dissipation produces a remainder $|\xi|^{2k-4}|\widehat{F_1}|^2$, losing the low-frequency weight needed for decay. We therefore retain the product estimate

$$|\xi|^{2k-2}|\widehat{F_1}|\,|\widehat{n^\varepsilon}| \lesssim |\xi|^{k-1}|\widehat{F_1}|\,G.$$

For the stress term involving $F_2$, controlling the stress factor directly through $G$ would introduce $\varepsilon^{-1/2}$, since the stress component of $G^2$ is weighted by $\varepsilon$. We therefore estimate this term using Young's inequality and the stress dissipation.

Together with the estimates for the $\pi_t$ terms, these two treatments bound the source by a small multiple of the dissipation plus $A_kG + B_k^2$, as shown in (3.22) and (3.23). After absorption, Lemma 3.3 yields

$$\partial_t G^2 + 2c\lambda_\xi G^2 \le C_0 A_k G + C_0 B_k^2, \qquad \lambda_\xi = \frac{|\xi|^2}{1+|\xi|^2},$$

where $A_k$ and $B_k$ are defined in (3.8), and $c, C_0 > 0$ are independent of $t$, $\xi$ and $\varepsilon$. Lemma A.2 then bounds $G$ by the damped initial-data term, a time convolution of $A_k$, and the square root of a time convolution of $B_k^2$.

The remaining difficulty is to close the nonlinear convolution estimates without assuming that the background is small. The source terms depend on the unknown perturbations, and their coupling with the background cannot be absorbed by smallness alone. To estimate these terms in terms of the same decay bounds that we seek to prove, we introduce the time-weighted functional

$$\begin{aligned} Y(t) = \sup_{\tau\in[0,t]}\Bigg\{ &\sum_{k=0}^{[\frac{d}{2}]+1} (1+\tau)^{\frac{d+2k+2}{4}}\|D^k(z^\varepsilon(\tau), n^\varepsilon(\tau))\| + (1+\tau)^{\frac{d+2[\frac{d}{2}]+4}{4}}\|D^{[\frac{d}{2}]+2}(z^\varepsilon(\tau), n^\varepsilon(\tau))\| \\ &+ (1+\tau)^{\frac{d+1}{2}}\|(z^\varepsilon(\tau), n^\varepsilon(\tau))\|_{L^\infty}\Bigg\}. \end{aligned}$$

The weighted derivative and $L^\infty$ bounds jointly control the nonlinear products. Frequency and time decompositions, together with $(\mathrm{H}_2)$ and interpolation using the uniform Sobolev bound, yield sublinear powers of $Y(t)$ from the background coupling terms. Young's inequality bounds these contributions by a constant plus an arbitrarily small multiple of $Y(t)$. The quadratic perturbation terms retain the small factor $\varepsilon^\alpha$. After absorption, Proposition 3.1 therefore gives

$$Y(t) \lesssim 1 + \varepsilon^\alpha Y(t)^2.$$

For each fixed $\alpha \in (0,1]$, Strauss' lemma [29] (Lemma A.1) gives $Y(t) \lesssim 1$ for sufficiently small $\varepsilon$. Finally, (1.6) transfers the momentum estimates to the velocity, and (2.4) gives the stress estimates after subtraction of the initial layer, completing the proof of Theorem 1.2.

The remainder of this paper is organized as follows. Section 2 proves the global well-posedness and convergence rate estimates in Theorem 1.1. Section 3 establishes the refined time decay estimates in Theorem 1.2. Appendix A contains two auxiliary lemmas used in Section 3.

**Notations.** We write $x = (x_1, \dots, x_d)^T \in \mathbb{R}^d$, $u = (u_1, \dots, u_d)^T$ and $\pi = (\pi_{ij})_{1\le i,j\le d}$. For vector and matrix fields, we use the conventions $(\nabla u)_{ij} := \partial_{x_j} u_i$, $(\operatorname{div}\pi)_i := \sum_{j=1}^d \partial_{x_j}\pi_{ij}$, $\operatorname{tr}(\pi) := \sum_{j=1}^d \pi_{jj}$ , and $\pi : \tau := \sum_{i,j=1}^d \pi_{ij}\tau_{ij}$. For a multi-index $\beta = (\beta_1, \dots, \beta_d)$, we write $|\beta| = \beta_1 + \cdots + \beta_d$ and $\partial_x^\beta = \partial_{x_1}^{\beta_1} \cdots \partial_{x_d}^{\beta_d}$. For each nonnegative integer $k$, $D^k f$ denotes the collection of $\partial_x^\beta f$ with $|\beta| = k$. For complex-valued $f(t,x)$, we denote its real part, imaginary part and conjugate by $\operatorname{Re}(f(t,x))$, $\operatorname{Im}(f(t,x))$ and $\overline{f(t,x)}$, respectively. Let i denote the imaginary unit. For $\xi = (\xi_1, \dots, \xi_d)^T \in \mathbb{R}^d$, the Fourier transform is

$$\widehat{f}(t,\xi) = \mathcal{F}(f)(t,\xi) := \int_{\mathbb{R}^d} e^{-\mathrm{i}x\cdot\xi} f(t,x)\,dx.$$

The Lebesgue space $L^p = L^p(\mathbb{R}^d)$, $1 \le p \le \infty$, has norm

$$\|f(t)\|_{L^p} := \begin{cases} \left(\displaystyle\int_{\mathbb{R}^d} |f(t,x)|^p\,dx\right)^{1/p}, & 1 \le p < \infty, \\ \operatorname*{ess\,sup}_{x\in\mathbb{R}^d} |f(t,x)|, & p = \infty. \end{cases}$$

In $L^2$, we use the inner product $\langle f(t), g(t)\rangle = \int_{\mathbb{R}^d} (f(t,x))^T \overline{g(t,x)} dx$ and abbreviate $\|f(t)\| := \|f(t)\|_{L^2}$. The norm $\|\widehat{f}(t)\|_{L^p_\xi}$ is taken with respect to the frequency variable $\xi$. For a nonnegative integer $s$, the Sobolev space $H^s = H^s(\mathbb{R}^d)$ has norm

$$\|f(t)\|_s^2 := \sum_{0\le k\le s} \|D^k f(t)\|^2.$$

For matrix fields, we use the conventions

$$\langle \pi(t), \tau(t)\rangle = \int_{\mathbb{R}^d} \pi(t,x) : \overline{\tau(t,x)} dx, \quad \|\pi(t)\|^2 = \langle \pi(t), \pi(t)\rangle, \quad \|\pi(t)\|_s^2 = \sum_{0\le k\le s} \|D^k\pi(t)\|^2.$$

Unless otherwise stated, the positive constants $C$, $C_j$ and $c_j$ ($j \in \mathbb{N}$) may depend on the fixed parameters, background bounds and uniform initial-data bounds, but are independent of $t$, $x$, $\xi$, $\varepsilon$ and the bootstrap bound $\delta$. In $C(\cdot)$ and $C_j(\cdot)$, the displayed arguments indicate additional dependencies. The constants $\epsilon > 0$ and $\epsilon_j > 0$ may be chosen sufficiently small, with the same convention on dependencies. We write $A \lesssim B$ if $A \le CB$ with $C$ as above, and $A \sim B$ if $A \lesssim B$ and $B \lesssim A$.

# 2 Proof of Theorem 1.1: Global Well-Posedness and Convergence Rate Estimates

In this section, we prove Theorem 1.1 by establishing uniform a priori estimates. We fix $\alpha \in (0,1]$, the coefficients, and the background solution. Dependence on these fixed quantities is suppressed below.

We first record the following local well-posedness result.

**Lemma 2.1.** *Under the assumptions of Theorem 1.1 on $(\rho(t,x), u(t,x))$, fix $\alpha \in (0,1]$ and let $s \in \mathbb{N}$ satisfy $\frac{d}{2}+1 < s \leq \frac{s_1-[d/2]-9}{2}$. Suppose that $(z_0^\varepsilon(x), w_0^\varepsilon(x), \chi_0^\varepsilon(x)) \in H^s$ and $\|(z_0^\varepsilon, w_0^\varepsilon, \sqrt{\varepsilon}\,\chi_0^\varepsilon)\|_s \leq M_0$, where $M_0$ is independent of $\varepsilon$. There exist $\varepsilon_{\rm loc} \in (0,1]$ and $T_0 > 0$, depending on $M_0$ and the fixed parameters and background bounds, such that, for every $0 < \varepsilon \leq \varepsilon_{\rm loc}$, the Cauchy problem (1.7), (1.8) admits a unique smooth solution $(z^\varepsilon(t,x), w^\varepsilon(t,x), \chi^\varepsilon(t,x)) \in C\left([0,T_0]; H^s(\mathbb{R}^d)\right) \cap C^1\left([0,T_0]; H^{s-1}(\mathbb{R}^d)\right)$ satisfying*

$$\|(z^\varepsilon(t), w^\varepsilon(t), \sqrt{\varepsilon}\chi^\varepsilon(t))\|_s \leq C_{\rm loc}\left(\|(z_0^\varepsilon, w_0^\varepsilon, \sqrt{\varepsilon}\chi_0^\varepsilon)\|_s + \varepsilon^{1-\alpha}\right), \quad 0 \leq t \leq T_0, \tag{2.1}$$

*where $C_{\rm loc}$ is a constant independent of $t, x$ and $\varepsilon$.*

To apply the local theory, we choose $\varepsilon_{\rm loc}$ sufficiently small so that Sobolev embedding ensures $\rho_0^\varepsilon(x) \geq \underline{\rho}/2$. Since (1.7) is symmetrized by $\mathrm{diag}(h'(\rho^\varepsilon), \rho^\varepsilon I_d)$, together with $\pi_t \in H^s$ under the assumed background regularity, this allows us to apply the standard local theory for symmetrizable hyperbolic systems. See [15] and [19, Section 2.1] for the corresponding proofs. Using the stress dissipation, we derive energy estimates in the weighted norm $\|(z^\varepsilon(t), w^\varepsilon(t), \sqrt{\varepsilon}\chi^\varepsilon(t))\|_s$ and obtain (2.1) with a lifespan $T_0$ independent of $\varepsilon$. We omit the details.

Assume that the above local solution $(z^\varepsilon(t,x), w^\varepsilon(t,x), \chi^\varepsilon(t,x)) \in C\left([0,T]; H^s(\mathbb{R}^d)\right) \cap C^1\left([0,T]; H^{s-1}(\mathbb{R}^d)\right)$ has been extended to time $T$. We impose the following a priori assumption for some $\delta > 0$:

$$\|(z^\varepsilon(t), w^\varepsilon(t), \sqrt{\varepsilon}\chi^\varepsilon(t))\|_s \leq \delta, \quad 0 \leq t \leq T. \tag{2.2}$$

The main a priori estimate of this section is the following.

**Proposition 2.2.** *Assume that the solution in Lemma 2.1 exists on $[0,T]$ and satisfies (2.2). There exists a sufficiently small $\varepsilon_0 \in (0, \varepsilon_{\rm loc}]$, depending only on $\delta$, such that, for $0 < \varepsilon \leq \varepsilon_0$ and all $t \in [0,T]$,*

$$\begin{aligned} &\|(z^\varepsilon(t), w^\varepsilon(t), \sqrt{\varepsilon}\chi^\varepsilon(t))\|_s^2 + \int_0^t \left(\|\nabla(z^\varepsilon(\tau), w^\varepsilon(\tau))\|_{s-1}^2 + \|\chi^\varepsilon(\tau)\|_s^2\right) d\tau \\ \lesssim &\|(z_0^\varepsilon, w_0^\varepsilon, \sqrt{\varepsilon}\chi_0^\varepsilon)\|_s^2 + \varepsilon^{2-2\alpha}. \end{aligned} \tag{2.3}$$

With the above results, we prove Theorem 1.1.

*Proof of Theorem 1.1.* By Proposition 2.2, choosing $\varepsilon$ sufficiently small closes the standard continuation argument, yielding global existence and (1.10). To obtain the remaining stress estimate, we integrate the Maxwell equation $(1.7)_3$ and find

$$\begin{aligned} \chi^\varepsilon(t,x) = &\, e^{-\frac{t}{\varepsilon}}\chi_0^\varepsilon(x) \\ &- \varepsilon^{-1}\int_0^t e^{-\frac{t-\tau}{\varepsilon}}\left(\nu\sigma(w^\varepsilon) + \lambda(\mathrm{div}\, w^\varepsilon)I_d + \varepsilon^{1-\alpha}\pi_t\right)(\tau,x)\, d\tau. \end{aligned} \tag{2.4}$$

Since $\varepsilon^{-1}\int_0^t e^{-(t-\tau)/\varepsilon}\,d\tau \le 1$, Minkowski's inequality, (2.3) and ($H_2$) yield (1.11). Together with (1.5), these estimates imply (1.12). □

To prove Proposition 2.2, we derive a series of energy estimates and recover the dissipation of the density and velocity perturbations.

Assume that

$$0<\varepsilon\le\hat{\varepsilon}_1:=\min\left\{1,\left(\frac{\underline{\rho}}{2C_s\delta}\right)^{1/\alpha},\left(\frac{\overline{\rho}}{C_s\delta}\right)^{1/\alpha},\left(\frac{\overline{u}}{C_s\delta}\right)^{1/\alpha}\right\},$$

where $C_s$ denotes the Sobolev embedding constant. By Sobolev embedding, (1.5), Assumption ($H_1$) and (2.2) give uniform bounds on $\rho^\varepsilon$ and $u^\varepsilon$. Together with the smoothness of $h(\rho)$, these imply

$$\begin{aligned}&\underline{\rho}/2\le\rho^\varepsilon(t,x)\le 2\overline{\rho},\quad |u^\varepsilon(t,x)|\le 2\overline{u},\\ &|h^{(k)}(\rho(t,x))|+|h^{(k)}(\rho^\varepsilon(t,x))|\lesssim 1,\quad 0\le k\le s,\ (t,x)\in[0,T]\times\mathbb{R}^d,\end{aligned}\tag{2.5}$$

where $h^{(k)}(\rho)=\frac{d^kh(\rho)}{d\rho^k}$.

We now establish the following $L^2$ energy estimate.

**Lemma 2.3.** *Assume that the solution obtained in Lemma 2.1 has been extended to $t=T$, and that $(z^\varepsilon(t,x),w^\varepsilon(t,x),\chi^\varepsilon(t,x))$ satisfies the a priori assumption (2.2). Then there exists a sufficiently small positive constant $\hat{\varepsilon}_1$, depending only on $\delta$, such that, whenever $0<\varepsilon\le\hat{\varepsilon}_1$,*

$$\|(z^\varepsilon(t),w^\varepsilon(t),\sqrt{\varepsilon}\chi^\varepsilon(t))\|^2+\int_0^t\|\chi^\varepsilon(\tau)\|^2d\tau\lesssim\|(z_0^\varepsilon,w_0^\varepsilon,\sqrt{\varepsilon}\chi_0^\varepsilon)\|^2+\varepsilon^{2-2\alpha},\quad 0\le t\le T.\tag{2.6}$$

*Proof.* Recall that $H'(\rho)=h(\rho)$ and $h'(\rho)=p'(\rho)/\rho$. We use the quadratic stress energy,

$$P_1(\pi)=\frac{1}{4\nu}\left(|\pi|^2+\omega(\operatorname{tr}\pi)^2\right),\qquad \omega=\frac{2\nu}{d^2\lambda}-\frac{1}{d}.$$

The inequality $(\operatorname{tr}\pi)^2\le d|\pi|^2$ yields

$$P_1(\pi)\ge\min\left\{\frac{1}{4\nu},\frac{1}{2d\lambda}\right\}|\pi|^2.\tag{2.7}$$

To derive the basic energy estimate, we define the relative entropy to the background by

$$\eta(t,x)=\frac12\rho^\varepsilon|u^\varepsilon-u|^2+H(\rho^\varepsilon)-H(\rho)-h(\rho)(\rho^\varepsilon-\rho)+\varepsilon P_1(\pi^\varepsilon-\pi).$$

In addition, we use the conservative variables $W=(\rho,m)^T=(\rho,\rho u)^T$ and $W^\varepsilon=(\rho^\varepsilon,m^\varepsilon)^T=(\rho^\varepsilon,\rho^\varepsilon u^\varepsilon)^T$, the Euler entropy $\eta_0(W)=|m|^2/(2\rho)+H(\rho)$, and the notation

$$F_j(W)=\begin{pmatrix}m_j\\ \frac{mm_j}{\rho}+p(\rho)e_j\end{pmatrix},\qquad A_j(W)=\nabla_W F_j(W),$$
$$A_j(W^\varepsilon,W)=\int_0^1 A_j(\theta W^\varepsilon+(1-\theta)W)\,d\theta.$$

Using (1.1), (1.3) and $\operatorname{div}(\pi^T u) = u \cdot \operatorname{div}\pi + \pi : \nabla u$, a direct calculation yields

$$\begin{aligned}
&\partial_t \eta + 2P_1(\pi^\varepsilon - \pi) + \operatorname{div}\mathcal{Q} \\
&= \sum_{j=1}^d W_{x_j}^T \nabla_W^2 \eta_0(W)\big(A_j(W) - A_j(W^\varepsilon, W)\big)(W^\varepsilon - W) + \frac{\rho^\varepsilon - \rho}{\rho}(\operatorname{div}\pi)\cdot(u^\varepsilon - u) \\
&\quad - \frac{\varepsilon}{2\nu}(\pi^\varepsilon - \pi) : \pi_t - \frac{\varepsilon\omega}{2\nu}\operatorname{tr}(\pi^\varepsilon - \pi)\operatorname{tr}(\pi_t).
\end{aligned} \tag{2.8}$$

Here

$$\begin{aligned}
\mathcal{Q}(t,x) = &\frac{1}{2}\rho^\varepsilon u^\varepsilon |u^\varepsilon - u|^2 + \left(H(\rho^\varepsilon) - H(\rho) - h(\rho)(\rho^\varepsilon - \rho)\right)u^\varepsilon \\
&+ \big(p(\rho^\varepsilon) - p(\rho)\big)(u^\varepsilon - u) + (\pi^\varepsilon - \pi)(u^\varepsilon - u).
\end{aligned}$$

Integrating (2.8) over $[0,t]\times\mathbb{R}^d$, noting that $A_j(W,W) = A_j(W)$, and applying Young's inequality, we obtain

$$\begin{aligned}
&\int_{\mathbb{R}^d} \eta(t,x)\,dx + \int_0^t \int_{\mathbb{R}^d} P_1(\pi^\varepsilon - \pi)(\tau,x)\,dxd\tau \\
\lesssim &\int_{\mathbb{R}^d} \eta(0,x)\,dx + \int_0^t \|\nabla(\rho(\tau),u(\tau))\|_{W^{1,\infty}} \|(\rho^\varepsilon(\tau) - \rho(\tau), u^\varepsilon(\tau) - u(\tau))\|^2\,d\tau \\
&+ \int_0^t \int_{\mathbb{R}^d} \varepsilon^2(|\pi_t(\tau,x)|^2 + |\operatorname{tr}(\pi_t(\tau,x))|^2)\,dxd\tau.
\end{aligned} \tag{2.9}$$

With (2.9) at hand, we observe that, for $0 < \varepsilon \le \hat{\varepsilon}_1$, (2.5) implies

$$\eta(t,x) \sim \varepsilon^{2\alpha}|(z^\varepsilon(t,x), w^\varepsilon(t,x))|^2 + \varepsilon^{1+2\alpha}|\chi^\varepsilon(t,x)|^2.$$

Using ($H_2$), (2.7) and Gronwall's inequality, we obtain (2.6). □

For the higher order energy estimate, we have the following lemma.

**Lemma 2.4.** *Under the assumptions of Lemma 2.3, for all $0 \le t \le T$, we have*

$$\begin{aligned}
&\|(z^\varepsilon(t), w^\varepsilon(t), \sqrt{\varepsilon}\chi^\varepsilon(t))\|_s^2 + \int_0^t \|\chi^\varepsilon(\tau)\|_s^2 d\tau \\
\le &C_0(\epsilon_1)\left(\|(z_0^\varepsilon, w_0^\varepsilon, \sqrt{\varepsilon}\chi_0^\varepsilon)\|_s^2 + \varepsilon^{2-2\alpha}\right) + (\epsilon_1 + C_1(\delta)\varepsilon^\alpha)\int_0^t \|D\left(z^\varepsilon(\tau), w^\varepsilon(\tau)\right)\|_{s-1}^2 d\tau,
\end{aligned} \tag{2.10}$$

*where $\epsilon_1 > 0$ can be chosen arbitrarily small.*

*Proof.* We first estimate $\chi^\varepsilon(t,x)$. For a multi-index $\beta \in \mathbb{N}^d$ with $1 \le |\beta| \le s$, applying $\partial_x^\beta$ to $(1.7)_3$ and taking the Frobenius product with the matrix $(2\nu)^{-1}\partial_x^\beta\chi^\varepsilon$, we obtain

$$\begin{aligned}
&\frac{\varepsilon}{4\nu}\frac{d}{dt}|\partial_x^\beta\chi^\varepsilon|^2 + \frac{1}{2\nu}|\partial_x^\beta\chi^\varepsilon|^2 + \partial_x^\beta\chi^\varepsilon : \nabla\partial_x^\beta w^\varepsilon \\
&\qquad + \left(\frac{\lambda}{2\nu} - \frac{1}{d}\right)(\operatorname{div}\partial_x^\beta w^\varepsilon)\cdot\operatorname{tr}(\partial_x^\beta\chi^\varepsilon) = -\frac{1}{2\nu}\varepsilon^{1-\alpha}\partial_x^\beta\pi_t : \partial_x^\beta\chi^\varepsilon.
\end{aligned} \tag{2.11}$$

On the other hand, applying $\partial_x^\beta$ to $(1.7)_3$, taking the trace, and then multiplying by $\operatorname{tr}(\partial_x^\beta \chi^\varepsilon)$, we use $\operatorname{tr}(\sigma(\partial_x^\beta w^\varepsilon)) = 0$ to obtain

$$\frac{\varepsilon}{2}\frac{d}{dt}|\operatorname{tr}(\partial_x^\beta \chi^\varepsilon)|^2 + |\operatorname{tr}(\partial_x^\beta \chi^\varepsilon)|^2 + \lambda d(\operatorname{div}\partial_x^\beta w^\varepsilon)\cdot \operatorname{tr}(\partial_x^\beta \chi^\varepsilon) = \varepsilon^{1-\alpha}\lambda d(\operatorname{div}\partial_x^\beta u)_t \cdot \operatorname{tr}(\partial_x^\beta \chi^\varepsilon). \tag{2.12}$$

Combining (2.11) and (2.12) and integrating over $\mathbb{R}^d$, we obtain

$$2\varepsilon\frac{d}{dt}\int_{\mathbb{R}^d} P_1(\partial_x^\beta \chi^\varepsilon)dx + 4\int_{\mathbb{R}^d} P_1(\partial_x^\beta \chi^\varepsilon)dx + 2\int_{\mathbb{R}^d} \partial_x^\beta \chi^\varepsilon : \nabla \partial_x^\beta w^\varepsilon dx = \int_{\mathbb{R}^d} \widetilde{F}_3 dx, \tag{2.13}$$

where

$$\widetilde{F}_3 = -\nu^{-1}\varepsilon^{1-\alpha}\partial_x^\beta \pi_t : \partial_x^\beta \chi^\varepsilon - 2\left(\frac{\lambda}{2\nu} - \frac{1}{d}\right)\varepsilon^{1-\alpha}(\operatorname{div}\partial_x^\beta u)_t \cdot \operatorname{tr}(\partial_x^\beta \chi^\varepsilon).$$

We next derive the corresponding estimate for $(z^\varepsilon(t,x), w^\varepsilon(t,x))$. To this end, let

$$U^\varepsilon(t,x) = \begin{pmatrix} z^\varepsilon(t,x) \\ w^\varepsilon(t,x) \end{pmatrix}, \quad B_j(\rho^\varepsilon, u^\varepsilon) = \begin{pmatrix} u_j + \varepsilon^\alpha w_j^\varepsilon & \rho^\varepsilon e_j^T \\ h'(\rho^\varepsilon)e_j & (u_j + \varepsilon^\alpha w_j^\varepsilon) I_d \end{pmatrix},$$

and

$$\widetilde{F} = \begin{pmatrix} \widetilde{F}_1 \\ \widetilde{F}_2 \end{pmatrix}, \quad S = \begin{pmatrix} 0 \\ -\frac{\operatorname{div}\chi^\varepsilon}{\rho^\varepsilon} \end{pmatrix}.$$

Then equations $(1.7)_{1,2}$ can be rewritten in the following matrix form:

$$\partial_t U^\varepsilon + \sum_{j=1}^d B_j(\rho^\varepsilon, u^\varepsilon)\partial_{x_j} U^\varepsilon = \widetilde{F} + S. \tag{2.14}$$

Observe that system (2.14) admits the symmetrizer

$$B_0(\rho^\varepsilon) = \operatorname{diag}(h'(\rho^\varepsilon), \rho^\varepsilon I_d).$$

Applying $\partial_x^\beta$ to (2.14), taking the $L^2$ inner product with $2B_0(\rho^\varepsilon)\partial_x^\beta U^\varepsilon$, and noting that $\widetilde{B}_j := B_0 B_j$ is symmetric, we apply the energy method for symmetric hyperbolic systems (see [19]) to obtain

$$\begin{aligned}
&\frac{d}{dt}\langle B_0(\rho^\varepsilon)\partial_x^\beta U^\varepsilon, \partial_x^\beta U^\varepsilon\rangle + 2\langle \partial_x^\beta w^\varepsilon, \operatorname{div}\partial_x^\beta \chi^\varepsilon\rangle \\
=&2\langle B_0(\rho^\varepsilon)\partial_x^\beta U^\varepsilon, \partial_x^\beta \widetilde{F}\rangle - 2\sum_{j=1}^d \langle B_0(\rho^\varepsilon)\partial_x^\beta U^\varepsilon, \partial_x^\beta(B_j\partial_{x_j}U^\varepsilon) - B_j\partial_x^\beta\partial_{x_j}U^\varepsilon\rangle \\
&- 2\left\langle \rho^\varepsilon\partial_x^\beta w^\varepsilon, \partial_x^\beta\left(\frac{\operatorname{div}\chi^\varepsilon}{\rho^\varepsilon}\right) - \frac{\partial_x^\beta \operatorname{div}\chi^\varepsilon}{\rho^\varepsilon}\right\rangle + \langle \operatorname{div}\vec{B}(U^\varepsilon)\partial_x^\beta U^\varepsilon, \partial_x^\beta U^\varepsilon\rangle,
\end{aligned} \tag{2.15}$$

where

$$\operatorname{div}\vec{B}(U^\varepsilon) = \operatorname{diag}(h''(\rho^\varepsilon), I_d)\partial_t \rho^\varepsilon + \sum_{j=1}^d \partial_{x_j}\widetilde{B}_j(U^\varepsilon).$$

Combining (2.13) and (2.15), we obtain

$$
\begin{aligned}
&\frac{d}{dt}\int_{\mathbb{R}^d} h'(\rho^\varepsilon)|\partial_x^\beta z^\varepsilon|^2 + \rho^\varepsilon|\partial_x^\beta w^\varepsilon|^2 + 2\varepsilon P_1(\partial_x^\beta \chi^\varepsilon)dx + 4\int_{\mathbb{R}^d} P_1(\partial_x^\beta \chi^\varepsilon)dx \\
=&2\langle h'(\rho^\varepsilon)\partial_x^\beta z^\varepsilon, \partial_x^\beta \widetilde{F}_1\rangle + 2\langle \rho^\varepsilon \partial_x^\beta w^\varepsilon, \partial_x^\beta \widetilde{F}_2\rangle - 2\sum_{j=1}^d K_j^\beta \\
&-2\left\langle \rho^\varepsilon \partial_x^\beta w^\varepsilon, \partial_x^\beta\left(\frac{\operatorname{div}\chi^\varepsilon}{\rho^\varepsilon}\right) - \frac{\operatorname{div}\partial_x^\beta \chi^\varepsilon}{\rho^\varepsilon}\right\rangle + \langle \operatorname{div}\vec{B}(U^\varepsilon)\partial_x^\beta U^\varepsilon, \partial_x^\beta U^\varepsilon\rangle + \int_{\mathbb{R}^d}\widetilde{F}_3 dx := \sum_{j=1}^6 I_j,
\end{aligned}
\tag{2.16}
$$

where

$$
\begin{aligned}
K_j^\beta =&\langle h'(\rho^\varepsilon)\partial_x^\beta z^\varepsilon, \partial_x^\beta(u_j^\varepsilon \partial_{x_j} z^\varepsilon) - u_j^\varepsilon \partial_x^\beta \partial_{x_j} z^\varepsilon\rangle + \langle h'(\rho^\varepsilon)\partial_x^\beta z^\varepsilon, \partial_x^\beta(\rho^\varepsilon \partial_{x_j} w_j^\varepsilon) - \rho^\varepsilon \partial_x^\beta \partial_{x_j} w_j^\varepsilon\rangle \\
&+\langle \rho^\varepsilon \partial_x^\beta w_j^\varepsilon, \partial_x^\beta(h'(\rho^\varepsilon)\partial_{x_j} z^\varepsilon) - h'(\rho^\varepsilon)\partial_x^\beta \partial_{x_j} z^\varepsilon\rangle + \langle \rho^\varepsilon \partial_x^\beta w^\varepsilon, \partial_x^\beta(u_j^\varepsilon \partial_{x_j} w^\varepsilon) - u_j^\varepsilon \partial_x^\beta \partial_{x_j} w^\varepsilon\rangle.
\end{aligned}
$$

We now estimate $I_1$-$I_6$. In view of $(1.1)_1$, (2.2), the Moser-type inequalities and Young's inequality, we obtain

$$
\begin{aligned}
I_1 &= -2\langle h'(\rho^\varepsilon)\partial_x^\beta z^\varepsilon, \partial_x^\beta(z^\varepsilon \operatorname{div} u) + \partial_x^\beta(w^\varepsilon \cdot \nabla \rho)\rangle \\
&\lesssim \left(\left\|D^2\left(\rho(t), u(t)\right)\right\|_{s-1}^2 + \|\nabla(\rho(t), u(t))\|_{L^\infty}^2\right)\|(z^\varepsilon(t), w^\varepsilon(t))\|_s^2 + \epsilon\|D(z^\varepsilon(t), w^\varepsilon(t))\|_{s-1}^2,
\end{aligned}
$$

$$
\begin{aligned}
I_2 =& 2\left\langle \rho^\varepsilon \partial_x^\beta w^\varepsilon, \partial_x^\beta\left(\frac{p'(\rho)z^\varepsilon \nabla\rho}{\rho\rho^\varepsilon}\right) + \partial_x^\beta\left(\frac{z^\varepsilon \operatorname{div}\pi}{\rho\rho^\varepsilon}\right) - \partial_x^\beta\left(\frac{z^\varepsilon \nabla\rho\int_0^1 p''(\rho+\theta\varepsilon^\alpha z^\varepsilon)d\theta}{\rho^\varepsilon}\right) - \partial_x^\beta(w^\varepsilon\cdot\nabla u)\right\rangle \\
\lesssim& \left(\epsilon + C(\delta)\varepsilon^\alpha\right)\|D(z^\varepsilon(t), w^\varepsilon(t))\|_{s-1}^2 + \left(\|D^2(\rho(t), u(t))\|_s^2 + \|\nabla(\rho(t), u(t))\|_{W^{1,\infty}}^2\right)\|(z^\varepsilon(t), w^\varepsilon(t))\|_s^2,
\end{aligned}
$$

$$
\begin{aligned}
I_3 \leq& 2\sum_{j=1}^d |K_j^\beta| \\
\lesssim& \left(\epsilon + C(\delta)\varepsilon^\alpha\right)\|D(z^\varepsilon(t), w^\varepsilon(t))\|_{s-1}^2 \\
&+\left(\|D(\rho(t), u(t))\|_{s-1}^2 + \|D(\rho(t), u(t))\|_{L^\infty}^2\right)\|(z^\varepsilon(t), w^\varepsilon(t))\|_s^2,
\end{aligned}
$$

$$
I_4 \lesssim \left(\epsilon + C(\delta)\varepsilon^\alpha\right)\|D(w^\varepsilon(t), \chi^\varepsilon(t))\|_{s-1}^2 + \left(\|D\rho(t)\|_{s-1}^2 + \|D\rho(t)\|_{L^\infty}^2\right)\|w^\varepsilon(t)\|_s^2,
$$

$$
I_5 \lesssim (\epsilon + C(\delta)\varepsilon^\alpha)\|D(z^\varepsilon(t), w^\varepsilon(t))\|_{s-1}^2 + \|\nabla(\rho(t), u(t))\|_{L^\infty}^2\|D(z^\varepsilon(t), w^\varepsilon(t))\|_{s-1}^2,
$$

and

$$
I_6 \lesssim \varepsilon^{2-2\alpha}\|D^2 u_t(t)\|_{s-1}^2 + \epsilon\|D\chi^\varepsilon(t)\|_{s-1}^2.
$$

Substituting the estimates for $I_1$-$I_6$ into (2.16), we have

$$
\frac{d}{dt}\int_{\mathbb{R}^d} h'(\rho^\varepsilon)|\partial_x^\beta z^\varepsilon|^2 + \rho^\varepsilon|\partial_x^\beta w^\varepsilon|^2 + 2\varepsilon P_1(\partial_x^\beta \chi^\varepsilon)dx + 4\int_{\mathbb{R}^d} P_1(\partial_x^\beta \chi^\varepsilon)dx
$$

$$
\begin{aligned}
&\lesssim \left(\|D(\rho(t),u(t))\|_{s+1}^2+\|D(\rho(t),u(t))\|_{W^{1,\infty}}^2\right)\|(z^\varepsilon(t),w^\varepsilon(t))\|_s^2\\
&\quad+(\epsilon+C(\delta)\varepsilon^\alpha)\|D(z^\varepsilon(t),w^\varepsilon(t),\chi^\varepsilon(t))\|_{s-1}^2+\varepsilon^{2-2\alpha}\|D^2u_t(t)\|_{s-1}^2.
\end{aligned}
$$

Summing the above inequality over $\beta$ with $1\le|\beta|\le s$ and integrating in time, we obtain

$$
\begin{aligned}
&\sum_{1\le|\beta|\le s}\left(\int_{\mathbb{R}^d}h'(\rho^\varepsilon)|\partial_x^\beta z^\varepsilon|^2+\rho^\varepsilon|\partial_x^\beta w^\varepsilon|^2+2\varepsilon P_1(\partial_x^\beta\chi^\varepsilon)dx+4\int_0^t\int_{\mathbb{R}^d}P_1(\partial_x^\beta\chi^\varepsilon)dxd\tau\right)\\
\lesssim&\sum_{1\le|\beta|\le s}\int_{\mathbb{R}^d}h'(\rho_0^\varepsilon)|\partial_x^\beta z_0^\varepsilon|^2+\rho_0^\varepsilon|\partial_x^\beta w_0^\varepsilon|^2+2\varepsilon P_1(\partial_x^\beta\chi_0^\varepsilon)dx\\
&+\int_0^t\left(\|D(\rho(\tau),u(\tau))\|_{s+1}^2+\|D(\rho(\tau),u(\tau))\|_{W^{1,\infty}}^2\right)\|(z^\varepsilon(\tau),w^\varepsilon(\tau))\|_s^2d\tau\\
&+(\epsilon+C(\delta)\varepsilon^\alpha)\int_0^t\|D(z^\varepsilon(\tau),w^\varepsilon(\tau),\chi^\varepsilon(\tau))\|_{s-1}^2\,d\tau+\varepsilon^{2-2\alpha}.
\end{aligned}
\tag{2.17}
$$

Combining (2.6), (2.7) and (2.17), then using ($H_2$) and Gronwall's inequality, we obtain (2.10). □

To close the a priori estimates, we need to derive dissipation estimates for $z^\varepsilon$ and $w^\varepsilon$. This is the main content of the following lemma.

**Lemma 2.5.** *Under the assumptions of Lemma 2.3, there exists a sufficiently small positive constant $\hat\varepsilon_2\le\hat\varepsilon_1$, depending only on $\delta$, such that, whenever $0<\varepsilon\le\hat\varepsilon_2$, for all $0\le t\le T$, we have*

$$
\begin{aligned}
&\int_0^t\|\nabla(z^\varepsilon(\tau),w^\varepsilon(\tau))\|_{s-1}^2d\tau\\
\le&C_2\Bigg(\|(z_0^\varepsilon,w_0^\varepsilon,\sqrt{\varepsilon}\chi_0^\varepsilon)\|_s^2+\|(z^\varepsilon(t),w^\varepsilon(t),\sqrt{\varepsilon}\chi^\varepsilon(t))\|_s^2+\varepsilon^{2-2\alpha}+\int_0^t\|\chi^\varepsilon(\tau)\|_s^2d\tau\\
&+\int_0^t\left(\|D(\rho(\tau),u(\tau))\|_{s+1}^2+\|(\rho(\tau)-1,u(\tau))\|_{W^{2,\infty}}^2\right)\|(z^\varepsilon(\tau),w^\varepsilon(\tau))\|_s^2d\tau\Bigg).
\end{aligned}
\tag{2.18}
$$

*Proof.* We first use $(1.7)_3$ to derive the dissipation estimate for $w^\varepsilon(t,x)$. For this purpose, set $\mu=\nu+\lambda-\frac{2\nu}{d}$ and let $\beta\in\mathbb{N}^d$ be a multi-index satisfying $|\beta|\le s-1$. Applying $\partial_x^\beta$ and the divergence operator to $(1.7)_3$, we obtain

$$
\varepsilon\partial_t\partial_x^\beta\operatorname{div}\chi^\varepsilon+\partial_x^\beta\operatorname{div}\chi^\varepsilon+\nu\Delta\partial_x^\beta w^\varepsilon+\mu\nabla(\operatorname{div}\partial_x^\beta w^\varepsilon)=-\varepsilon^{1-\alpha}\partial_x^\beta\operatorname{div}\pi_t. \tag{2.19}
$$

Taking the $L^2$ inner product of (2.19) with $\partial_x^\beta w^\varepsilon$ and integrating by parts, we obtain

$$
\begin{aligned}
&\nu\|\nabla\partial_x^\beta w^\varepsilon(t)\|^2+\mu\|\operatorname{div}(\partial_x^\beta w^\varepsilon)(t)\|^2\\
=&-\varepsilon^{1-\alpha}\langle\partial_x^\beta\pi_t,\nabla\partial_x^\beta w^\varepsilon\rangle-\langle\partial_x^\beta\chi^\varepsilon,\nabla\partial_x^\beta w^\varepsilon\rangle+\varepsilon\langle\partial_t(\partial_x^\beta\operatorname{div}\chi^\varepsilon),\partial_x^\beta w^\varepsilon\rangle\\
=&I_7+I_8+I_9.
\end{aligned}
\tag{2.20}
$$

We next estimate $I_7$-$I_9$ successively. Let $4\kappa=\min\{\nu,\lambda\}$. By Young's inequality,

$$
I_7\le\kappa\|\nabla(\partial_x^\beta w^\varepsilon)(t)\|^2+C\varepsilon^{2-2\alpha}\|\nabla u_t(t)\|_{s-1}^2.
$$

Moreover,

$$I_8 \le \kappa\|\nabla(\partial_x^\beta w^\varepsilon)(t)\|^2 + C\|\chi^\varepsilon(t)\|_{s-1}^2.$$

Using integration by parts, Young's inequality and $(1.7)_2$, we obtain

$$\begin{aligned} I_9 =&\varepsilon\frac{d}{dt}\langle\partial_x^\beta \operatorname{div}\chi^\varepsilon, \partial_x^\beta w^\varepsilon\rangle - \varepsilon\langle\partial_x^\beta \operatorname{div}\chi^\varepsilon, \partial_x^\beta w_t^\varepsilon\rangle \\ =&\varepsilon\frac{d}{dt}\langle\partial_x^\beta \operatorname{div}\chi^\varepsilon, \partial_x^\beta w^\varepsilon\rangle + \varepsilon\left\langle\partial_x^\beta \operatorname{div}\chi^\varepsilon, \partial_x^\beta\left\{u^\varepsilon\cdot\nabla w^\varepsilon + h'(\rho^\varepsilon)\nabla z^\varepsilon + \frac{\operatorname{div}\chi^\varepsilon}{\rho^\varepsilon}\right\}\right\rangle \\ &- \varepsilon\left\langle\partial_x^\beta \operatorname{div}\chi^\varepsilon, \partial_x^\beta \widetilde{F}_2\right\rangle \\ =:&\varepsilon\frac{d}{dt}\langle\partial_x^\beta \operatorname{div}\chi^\varepsilon, \partial_x^\beta w^\varepsilon\rangle + I_{91} + I_{92}. \end{aligned} \tag{2.21}$$

The Moser-type inequalities and Young's inequality yield

$$I_{91} \lesssim C(\delta)\varepsilon\|D(z^\varepsilon(t), w^\varepsilon(t))\|_{s-1}^2 + C(\delta)\varepsilon\|\chi^\varepsilon(t)\|_s^2,$$

and

$$I_{92} \lesssim \varepsilon\left(\|D(\rho(t), u(t))\|_s^2 + \|\nabla(\rho(t), u(t))\|_{W^{1,\infty}}^2\right)\|(z^\varepsilon(t), w^\varepsilon(t))\|_s^2 + C(\delta)\varepsilon\|\chi^\varepsilon(t)\|_s^2.$$

Substituting the estimates for $I_{91}$ and $I_{92}$ into (2.21), we obtain

$$\begin{aligned} I_9 \lesssim&\varepsilon\frac{d}{dt}\left\langle\partial_x^\beta \operatorname{div}\chi^\varepsilon, \partial_x^\beta w^\varepsilon\right\rangle + C(\delta)\varepsilon\left(\|\chi^\varepsilon(t)\|_s^2 + \|D(z^\varepsilon(t), w^\varepsilon(t))\|_{s-1}^2\right) \\ &+ \varepsilon\left(\|D(\rho(t), u(t))\|_s^2 + \|\nabla(\rho(t), u(t))\|_{W^{1,\infty}}^2\right)\|(z^\varepsilon(t), w^\varepsilon(t))\|_s^2. \end{aligned}$$

Substituting the estimates for $I_7$-$I_9$ into (2.20), we have

$$\begin{aligned} \kappa\|\nabla\partial_x^\beta w^\varepsilon(t)\|^2 \lesssim&\varepsilon\frac{d}{dt}\langle\partial_x^\beta \operatorname{div}\chi^\varepsilon, \partial_x^\beta w^\varepsilon\rangle + C(\delta)\varepsilon\|D(z^\varepsilon(t), w^\varepsilon(t))\|_{s-1}^2 + (C(\delta)\varepsilon + C)\|\chi^\varepsilon(t)\|_s^2 \\ &+ \varepsilon^{2-2\alpha}\|\nabla u_t(t)\|_{s-1}^2 + \varepsilon\left(\|D(\rho(t), u(t))\|_s^2 + \|\nabla(\rho(t), u(t))\|_{W^{1,\infty}}^2\right)\|(z^\varepsilon(t), w^\varepsilon(t))\|_s^2. \end{aligned}$$

Summing the above inequality over $\beta$ with $|\beta| \le s-1$ and integrating in time, we obtain

$$\begin{aligned} \int_0^t \|\nabla w^\varepsilon(\tau)\|_{s-1}^2 d\tau \le& C_3\Bigg(\|(w_0^\varepsilon, \sqrt{\varepsilon}\chi_0^\varepsilon)\|_s^2 + \varepsilon^{2-2\alpha} + \|(w^\varepsilon(t), \sqrt{\varepsilon}\chi^\varepsilon(t))\|_s^2 \\ &+ C_2(\delta)\varepsilon\int_0^t \|D(z^\varepsilon(\tau), w^\varepsilon(\tau))\|_{s-1}^2 + \|\chi^\varepsilon(\tau)\|_s^2 d\tau + \int_0^t \|\chi^\varepsilon(\tau)\|_s^2 d\tau \\ &+ \varepsilon\int_0^t \left(\|D(\rho(\tau), u(\tau))\|_s^2 + \|\nabla(\rho(\tau), u(\tau))\|_{W^{1,\infty}}^2\right)\|(z^\varepsilon(\tau), w^\varepsilon(\tau))\|_s^2 d\tau\Bigg). \end{aligned} \tag{2.22}$$

We next use the velocity equation $(1.7)_2$ to recover the dissipation of $z^\varepsilon$. Applying $\partial_x^\beta$ to this equation and taking the $L^2$ inner product with $\partial_x^\beta(h'(\rho^\varepsilon)\nabla z^\varepsilon)$, we get

$$\begin{aligned} \|\partial_x^\beta(h'(\rho^\varepsilon)\nabla z^\varepsilon)(t)\|^2 =& -\langle\partial_t\partial_x^\beta w^\varepsilon, \partial_x^\beta(h'(\rho^\varepsilon)\nabla z^\varepsilon)\rangle \\ &- \left\langle\partial_x^\beta\left\{u^\varepsilon\cdot\nabla w^\varepsilon + \frac{\operatorname{div}\chi^\varepsilon}{\rho^\varepsilon} - \widetilde{F}_2\right\}, \partial_x^\beta(h'(\rho^\varepsilon)\nabla z^\varepsilon)\right\rangle \\ =:& I_{10} + I_{11}. \end{aligned} \tag{2.23}$$

We estimate $I_{10}$ and $I_{11}$ in turn. Observe that

$$\partial_t(h'(\rho^\varepsilon)\nabla z^\varepsilon) = \partial_t h'(\rho^\varepsilon)\nabla z^\varepsilon + h'(\rho^\varepsilon)\nabla\partial_t z^\varepsilon = \partial_t h'(\rho^\varepsilon)\nabla z^\varepsilon + \nabla(h'(\rho^\varepsilon)\partial_t z^\varepsilon) - \nabla h'(\rho^\varepsilon)\partial_t z^\varepsilon.$$

We have

$$\begin{aligned} I_{10} = & -\frac{d}{dt}\langle\partial_x^\beta w^\varepsilon, \partial_x^\beta(h'(\rho^\varepsilon)\nabla z^\varepsilon)\rangle + \langle\partial_x^\beta w^\varepsilon, \partial_x^\beta(\partial_t h'(\rho^\varepsilon)\nabla z^\varepsilon)\rangle \\ & - \langle\partial_x^\beta \operatorname{div} w^\varepsilon, \partial_x^\beta(h'(\rho^\varepsilon)\partial_t z^\varepsilon)\rangle - \langle\partial_x^\beta w^\varepsilon, \partial_x^\beta(\nabla h'(\rho^\varepsilon)\partial_t z^\varepsilon)\rangle \\ =: & -\frac{d}{dt}\langle\partial_x^\beta w^\varepsilon, \partial_x^\beta(h'(\rho^\varepsilon)\nabla z^\varepsilon)\rangle + I_{101} + I_{102} + I_{103}. \end{aligned} \tag{2.24}$$

Owing to $(1.3)_1$, $(1.7)_1$, the Moser-type inequalities and Young's inequality, we obtain

$$\begin{aligned} I_{101} = & \langle\partial_x^\beta w^\varepsilon, \partial_x^\beta(h''(\rho^\varepsilon)\partial_t\rho\nabla z^\varepsilon)\rangle + \varepsilon^\alpha\langle\partial_x^\beta w^\varepsilon, \partial_x^\beta(h''(\rho^\varepsilon)\partial_t z^\varepsilon\nabla z^\varepsilon)\rangle \\ = & -\langle\partial_x^\beta w^\varepsilon, \partial_x^\beta(h''(\rho^\varepsilon)\operatorname{div}(\rho u)\nabla z^\varepsilon)\rangle \\ & - \varepsilon^\alpha\langle\partial_x^\beta w^\varepsilon, \partial_x^\beta(h''(\rho^\varepsilon)\{u^\varepsilon\cdot\nabla z^\varepsilon + \rho^\varepsilon\operatorname{div} w^\varepsilon + z^\varepsilon\operatorname{div} u + w^\varepsilon\cdot\nabla\rho\}\nabla z^\varepsilon)\rangle \\ \lesssim & \left(\|D(\rho(t),u(t))\|_{s-1}^2 + \|\nabla(\rho(t),u(t))\|_{L^\infty}^2\right)\|w^\varepsilon(t)\|_{s-1}^2 + (\epsilon + C(\delta)\varepsilon^\alpha)\|D\left(z^\varepsilon(t),w^\varepsilon(t)\right)\|_{s-1}^2, \end{aligned}$$

$$\begin{aligned} I_{102} = & \langle\partial_x^\beta \operatorname{div} w^\varepsilon, \partial_x^\beta(h'(\rho^\varepsilon)\{u^\varepsilon\cdot\nabla z^\varepsilon + (\rho^\varepsilon - 1)\operatorname{div} w^\varepsilon + \operatorname{div} w^\varepsilon + z^\varepsilon\operatorname{div} u + w^\varepsilon\cdot\nabla\rho\})\rangle \\ \lesssim & \|Dw^\varepsilon(t)\|_{s-1}^2 + (\epsilon + C(\delta)\varepsilon^\alpha)\|D(z^\varepsilon(t),w^\varepsilon(t))\|_{s-1}^2 \\ & + \left(\|D(\rho(t),u(t))\|_{s-1}^2 + \|(\rho(t)-1,u(t))\|_{W^{1,\infty}}^2\right)\|(z^\varepsilon(t),w^\varepsilon(t))\|_s^2, \end{aligned}$$

$$\begin{aligned} I_{103} = & -\langle\partial_x^\beta w^\varepsilon, \partial_x^\beta(h''(\rho^\varepsilon)\nabla\rho\partial_t z^\varepsilon)\rangle - \varepsilon^\alpha\langle\partial_x^\beta w^\varepsilon, \partial_x^\beta(h''(\rho^\varepsilon)\nabla z^\varepsilon\partial_t z^\varepsilon)\rangle \\ = & \langle\partial_x^\beta w^\varepsilon, \partial_x^\beta(h''(\rho^\varepsilon)\nabla\rho\{u^\varepsilon\cdot\nabla z^\varepsilon + (\rho^\varepsilon-1)\operatorname{div} w^\varepsilon + \operatorname{div} w^\varepsilon + z^\varepsilon\operatorname{div} u + w^\varepsilon\cdot\nabla\rho\})\rangle \\ & + \varepsilon^\alpha\langle\partial_x^\beta w^\varepsilon, \partial_x^\beta(h''(\rho^\varepsilon)\nabla z^\varepsilon\{u^\varepsilon\cdot\nabla z^\varepsilon + (\rho^\varepsilon-1)\operatorname{div} w^\varepsilon + \operatorname{div} w^\varepsilon + z^\varepsilon\operatorname{div} u + w^\varepsilon\cdot\nabla\rho\})\rangle \\ \lesssim & \|Dw^\varepsilon(t)\|_{s-1}^2 + (\epsilon + C(\delta)\varepsilon^\alpha)\|D(z^\varepsilon(t),w^\varepsilon(t))\|_{s-1}^2 \\ & + \left(\|D(\rho(t),u(t))\|_{s-1}^2 + \|(\rho(t)-1,u(t))\|_{W^{1,\infty}}^2\right)\|(z^\varepsilon(t),w^\varepsilon(t))\|_s^2. \end{aligned}$$

Substituting the estimates for $I_{101}$-$I_{103}$ into (2.24), we obtain

$$\begin{aligned} I_{10} \lesssim & -\frac{d}{dt}\langle\partial_x^\beta w^\varepsilon, \partial_x^\beta(h'(\rho^\varepsilon)\nabla z^\varepsilon)\rangle + \|Dw^\varepsilon(t)\|_{s-1}^2 + (\epsilon + C(\delta)\varepsilon^\alpha)\|D(z^\varepsilon(t),w^\varepsilon(t))\|_{s-1}^2 \\ & + \left(\|D(\rho(t),u(t))\|_{s-1}^2 + \|(\rho(t)-1,u(t))\|_{W^{1,\infty}}^2\right)\|(z^\varepsilon(t),w^\varepsilon(t))\|_s^2. \end{aligned}$$

The Moser-type inequalities and Young's inequality give

$$\begin{aligned} I_{11} \lesssim & \|D\chi^\varepsilon(t)\|_{s-1}^2 + (\epsilon + C(\delta)\varepsilon^\alpha)\|D(z^\varepsilon(t),w^\varepsilon(t),\chi^\varepsilon(t))\|_{s-1}^2 \\ & + \left(\|D(\rho(t),u(t))\|_{s+1}^2 + \|D\rho(t)\|_{L^\infty}^2 + \|u(t)\|_{W^{2,\infty}}^2\right)\|(z^\varepsilon(t),w^\varepsilon(t))\|_s^2. \end{aligned}$$

Substituting the estimates for $I_{10}$ and $I_{11}$ into (2.23), we obtain

$$\begin{aligned} & \frac{d}{dt}\langle\partial_x^\beta w^\varepsilon, \partial_x^\beta(h'(\rho^\varepsilon)\nabla z^\varepsilon)\rangle + \|\partial_x^\beta(h'(\rho^\varepsilon)\nabla z^\varepsilon)\|^2 \\ \lesssim & \left(\|D(\rho(t),u(t))\|_{s+1}^2 + \|(\rho(t)-1,u(t))\|_{W^{2,\infty}}^2\right)\|(z^\varepsilon(t),w^\varepsilon(t))\|_s^2 \\ & + \|\chi^\varepsilon(t)\|_s^2 + \|Dw^\varepsilon(t)\|_{s-1}^2 + (\epsilon + C(\delta)\varepsilon^\alpha)\|D(z^\varepsilon(t),w^\varepsilon(t),\chi^\varepsilon(t))\|_{s-1}^2. \end{aligned}$$

Summing the above inequality over $\beta$ with $|\beta| \le s-1$ and integrating in time, and noting that

$$\sum_{|\beta|\le s-1}\left(\langle \partial_x^\beta w^\varepsilon, \partial_x^\beta (h'(\rho^\varepsilon)\nabla z^\varepsilon)\rangle + \int_0^t \|\partial_x^\beta (h'(\rho^\varepsilon)\nabla z^\varepsilon)(\tau)\|^2 d\tau\right)$$
$$\gtrsim (1 - C_3(\delta)\varepsilon^\alpha)\int_0^t \|\nabla z^\varepsilon(\tau)\|_{s-1}^2 d\tau - (1 + C(\delta)\varepsilon^\alpha)\left(\|(z_0^\varepsilon, w_0^\varepsilon)\|_s^2 + \|(z^\varepsilon(t), w^\varepsilon(t))\|_s^2\right),$$

we have, whenever $0 < \varepsilon \le \tilde\varepsilon_2 := (2C_3(\delta))^{-\frac{1}{\alpha}}$,

$$\begin{aligned}
&\int_0^t \|\nabla z^\varepsilon(\tau)\|_{s-1}^2 d\tau \\
\le& C_4\Bigg((1 + C_4(\delta)\varepsilon^\alpha)\left(\|(z_0^\varepsilon, w_0^\varepsilon, \sqrt{\varepsilon}\chi_0^\varepsilon)\|_s^2 + \|(z^\varepsilon(t), w^\varepsilon(t), \sqrt{\varepsilon}\chi^\varepsilon(t))\|_s^2\right) \\
&+ \int_0^t \|\nabla w^\varepsilon(\tau)\|_{s-1}^2 d\tau + \int_0^t \left(\|D(\rho(\tau), u(\tau))\|_{s+1}^2 + \|(\rho(\tau)-1, u(\tau))\|_{W^{2,\infty}}^2\right)\|(z^\varepsilon(\tau), w^\varepsilon(\tau))\|_s^2 d\tau \\
&+ \int_0^t \|\chi^\varepsilon(\tau)\|_s^2 d\tau + (\epsilon_2 + C_5(\delta)\varepsilon^\alpha)\int_0^t \|D(z^\varepsilon(\tau), w^\varepsilon(\tau), \chi^\varepsilon(\tau))\|_{s-1}^2 d\tau\Bigg).
\end{aligned} \tag{2.25}$$

We combine (2.22) and (2.25). For a sufficiently small positive constant $\epsilon_3$ independent of $t, x, \varepsilon$, and $\delta$, we choose a positive constant $\tilde\varepsilon_3$ satisfying

$$\begin{aligned}
C_3C_2(\delta)\tilde\varepsilon_3 + C_4(\epsilon_3 + \epsilon_2\epsilon_3 + C_5(\delta)\tilde\varepsilon_3^\alpha\epsilon_3) &\le \frac{1}{2}, \\
C_3C_2(\delta)\tilde\varepsilon_3 + C_4\epsilon_3\left(\epsilon_2 + C_5(\delta)\tilde\varepsilon_3^\alpha\right) &\le \frac{\epsilon_3}{2}, \\
C_3 + C_3C_2(\delta)\tilde\varepsilon_3 + C_4C_4(\delta)\tilde\varepsilon_3^\alpha\epsilon_3 + 1 &\le C_2.
\end{aligned}$$

When $0 < \hat\varepsilon_2 < \min\{\hat\varepsilon_1, \tilde\varepsilon_2, \tilde\varepsilon_3\}$, we obtain (2.18). □

To conclude this section, we note that adding a sufficiently small multiple of (2.18) to (2.10) allows us to absorb the energy and dissipation terms for sufficiently small $\epsilon_1$ and $\varepsilon$. Assumption ($H_2$) ensures the time integrability of the remaining background coefficient. Gronwall's inequality then yields (2.3), which completes the proof of Proposition 2.2.

# 3 Proof of Theorem 1.2: Time Decay Estimates

In this section, we prove Theorem 1.2. We introduce the time-weighted functional

$$\begin{aligned}
Y(t) = \sup_{\tau\in[0,t]}\Bigg\{&\sum_{k=0}^{[\frac{d}{2}]+1}(1+\tau)^{\frac{d+2k+2}{4}}\|D^k(z^\varepsilon(\tau), n^\varepsilon(\tau))\| + (1+\tau)^{\frac{d+2[\frac{d}{2}]+4}{4}}\|D^{[\frac{d}{2}]+2}(z^\varepsilon(\tau), n^\varepsilon(\tau))\| \\
&+ (1+\tau)^{\frac{d+1}{2}}\|(z^\varepsilon(\tau), n^\varepsilon(\tau))\|_{L^\infty}\Bigg\}.
\end{aligned}$$

The main estimate of this section is the following.

**Proposition 3.1.** *Under the assumptions of Theorem 1.2, for each fixed $\alpha \in (0,1]$ and all sufficiently small $\varepsilon > 0$, the global solution constructed in Theorem 1.1 satisfies, for every $T > 0$,*

$$Y(t) \lesssim 1 + \varepsilon^{\alpha} Y(t)^2, \qquad 0 \le t \le T. \tag{3.1}$$

*Consequently,*

$$Y(t) \lesssim 1, \qquad t \ge 0.$$

*Proof of Theorem 1.2.* Proposition 3.1 gives $Y(t) \lesssim 1$. Together with the uniform Sobolev and density bounds from Theorem 1.1, relation (1.6), the Moser-type inequalities and ($H_2$) yield (1.14). The Maxwell representation (2.4), together with the background estimates for $\pi_t$ obtained from the Navier-Stokes equations, gives (1.15). □

We now establish the nonlinear estimate (3.1), beginning with a bound for the background forcing $\pi_t$.

**Lemma 3.2.** *Let $d \ge 2$ and let $(\rho(t,x), u(t,x))$ be the global smooth solution of the Navier-Stokes equations* (1.3) *satisfying* (H1) *and* (H2). *Then*

$$\sup_{t \ge 0} \big\| (\widehat{\rho(t)-1}, \widehat{u}(t)) \big\|_{L^\infty_\xi} \le C.$$

*Proof.* We rewrite the Navier-Stokes equations (1.3) as

$$\begin{cases} (\rho-1)_t + \nabla \cdot u = f, \\ u_t + h'(1)\nabla(\rho-1) - \nu\Delta u - \mu\nabla(\nabla\cdot u) = g, \end{cases} \tag{3.2}$$

where

$$\begin{aligned} f(t,x) =& -u\cdot\nabla\rho - (\rho-1)\nabla\cdot u, \\ g(t,x) =& -u\cdot\nabla u - \left(\frac{p'(\rho)}{\rho} - h'(1)\right)\nabla\rho - \frac{\rho-1}{\rho}\left[\nu\Delta u + \mu\nabla(\nabla\cdot u)\right]. \end{aligned}$$

By ($H_1$), the smoothness of $p(\rho)$ and Hölder's inequality, we have

$$\begin{aligned} \|(f(t),g(t))\|_{L^1} \lesssim& \|(\rho(t)-1,u(t))\|_{L^2}\|\nabla(\rho(t)-1,u(t))\|_{L^2} \\ &+ \|\rho(t)-1\|_{L^2}\|D^2u(t)\|_{L^2}. \end{aligned}$$

It follows from ($H_2$) that

$$\begin{aligned} \|(f(t),g(t))\|_{L^1} &\lesssim (1+t)^{-\frac{d+1}{2}} + (1+t)^{-\frac{d+2}{2}} \\ &\lesssim (1+t)^{-\frac{d+1}{2}}. \end{aligned}$$

Since $d \ge 2$, this yields

$$\int_0^\infty \|(f(\tau),g(\tau))\|_{L^1}\, d\tau \lesssim 1. \tag{3.3}$$

Taking the Fourier transform of (3.2), multiplying the resulting equations by $h'(1)\overline{\widehat{\rho-1}}$ and $\overline{\widehat{u}}$, respectively, and taking the real part of their sum, we obtain

$$\frac{1}{2}\frac{d}{dt}\left(h'(1)|\widehat{\rho-1}|^2 + |\widehat{u}|^2\right) + \nu|\xi|^2|\widehat{u}|^2 + \mu|\xi\cdot\widehat{u}|^2$$

$$= h'(1)\,\mathrm{Re}\left(\widehat{f}\,\overline{\widehat{\rho-1}}\right) + \mathrm{Re}\left(\widehat{g}\cdot\overline{\widehat{u}}\right).$$

Dropping the nonnegative dissipation terms and applying the Cauchy–Schwarz inequality, we deduce

$$\left|(\widehat{\rho-1},\widehat{u})\right| \lesssim \left|(\widehat{\rho_0-1},\widehat{u}_0)\right| + \int_0^t \left|(\widehat{f},\widehat{g})\right| d\tau.$$

Using (3.3) and the initial $L^1$ regularity in (H2), we obtain

$$\sup_{t\geq 0}\left\|(\widehat{\rho-1},\widehat{u})\right\|_{L^\infty_\xi} \lesssim \|(\rho_0-1,u_0)\|_{L^1} + \int_0^\infty \|(f(\tau),g(\tau))\|_{L^1}\,d\tau \lesssim 1.$$

This completes the proof. □

To derive the decay estimates (1.14) and (1.15), we need to preserve the conservative structure of the linearized equations. To this end, in view of (1.6), we find that $(z^\varepsilon(t,x),n^\varepsilon(t,x),\chi^\varepsilon(t,x))$ satisfies the following Cauchy problem:

$$\begin{cases} z^\varepsilon_t + \mathrm{div}\, n^\varepsilon = 0,\\ n^\varepsilon_t + h'(1)\nabla z^\varepsilon + \mathrm{div}\,\chi^\varepsilon = F_1,\\ \varepsilon\chi^\varepsilon_t + \chi^\varepsilon + \nu\sigma(n^\varepsilon) + \lambda(\mathrm{div}\, n^\varepsilon)I_d = -\varepsilon^{1-\alpha}\pi_t + F_2, \end{cases} \tag{3.4}$$

with the initial data

$$\begin{aligned} &(z^\varepsilon(t,x),n^\varepsilon(t,x),\chi^\varepsilon(t,x))\,|_{t=0}\\ =&(z^\varepsilon_0(x),n^\varepsilon_0(x),\chi^\varepsilon_0(x))\\ :=&\varepsilon^{-\alpha}(\rho^\varepsilon_0(x)-\rho_0(x),m^\varepsilon_0(x)-m_0(x),\pi^\varepsilon_0(x)-\pi_0(x))\\ \to&(0,0,0),\quad |x|\to+\infty, \end{aligned}$$

where

$$\begin{aligned} F_1 =& -\nabla\left((p'(\rho^\varepsilon,\rho)-p'(1))z^\varepsilon\right) - \mathrm{div}\left(\frac{m\otimes n^\varepsilon}{\rho^\varepsilon} + \frac{n^\varepsilon\otimes m}{\rho^\varepsilon} + \varepsilon^\alpha\frac{n^\varepsilon\otimes n^\varepsilon}{\rho^\varepsilon} - \frac{z^\varepsilon m\otimes m}{\rho\rho^\varepsilon}\right),\\ F_2 =& \nu\sigma\left(\frac{(\rho^\varepsilon-1)n^\varepsilon}{\rho^\varepsilon} + \frac{z^\varepsilon m}{\rho\rho^\varepsilon}\right) + \lambda\,\mathrm{div}\left(\frac{(\rho^\varepsilon-1)n^\varepsilon}{\rho^\varepsilon} + \frac{z^\varepsilon m}{\rho\rho^\varepsilon}\right)I_d,\\ p'(\rho^\varepsilon,\rho) =& \int_0^1 p'(\theta\rho^\varepsilon + (1-\theta)\rho)d\theta. \end{aligned}$$

Taking the Fourier transform of (3.4) with respect to $x$, we obtain

$$\begin{cases} \widehat{z^\varepsilon_t} + \mathrm{i}\xi\cdot\widehat{n^\varepsilon} = 0,\\ \widehat{n^\varepsilon_t} + h'(1)\mathrm{i}\xi\widehat{z^\varepsilon} + \mathrm{i}\widehat{\chi^\varepsilon}\xi = \widehat{F_1},\\ \varepsilon\widehat{\chi^\varepsilon_t} + \widehat{\chi^\varepsilon} + \nu\widehat{\sigma(n^\varepsilon)} + \lambda(\mathrm{i}\xi\cdot\widehat{n^\varepsilon})I_d = -\varepsilon^{1-\alpha}\widehat{\pi}_t + \widehat{F_2}, \end{cases} \tag{3.5}$$

where $\widehat{\sigma(n^\varepsilon)} = \mathrm{i}\widehat{n^\varepsilon}\xi^T + (\mathrm{i}\widehat{n^\varepsilon}\xi^T)^T - \frac{2}{d}(\mathrm{i}\xi\cdot\widehat{n^\varepsilon})I_d$.

The following lemma plays a crucial role in the proof of the time decay estimates in Theorem 1.2.

**Lemma 3.3.** *For each integer $k \geq 1$ and all sufficiently small $0 < \varepsilon \leq 1$, there exists a nonnegative function $G(t,\xi)$ satisfying*

$$G(t,\xi)^2 \sim |\xi|^{2k-2}(1+|\xi|^2)\left(|\widehat{z^\varepsilon}(t,\xi)|^2 + |\widehat{n^\varepsilon}(t,\xi)|^2\right) + \varepsilon|\xi|^{2k}|\widehat{\chi^\varepsilon}(t,\xi)|^2. \tag{3.6}$$

*Moreover, setting $\lambda_\xi = \frac{|\xi|^2}{1+|\xi|^2}$, we obtain constants $c, C_0 > 0$ such that*

$$\partial_t G(t,\xi)^2 + 2c\lambda_\xi G(t,\xi)^2 \leq C_0 A_k(t,\xi)G(t,\xi) + C_0 B_k(t,\xi)^2, \tag{3.7}$$

*where*

$$\begin{aligned} A_k(t,\xi) &= \left(|\xi|^{k-1} + |\xi|^k\right)|\widehat{F}_1(t,\xi)| + \varepsilon^{1-\alpha}|\xi|^k|\widehat{\pi}_t(t,\xi)|, \\ B_k(t,\xi)^2 &= \left(|\xi|^{2k-2} + |\xi|^{2k}\right)|\widehat{F}_2(t,\xi)|^2 + \varepsilon^{2-2\alpha}|\xi|^{2k}|\widehat{\pi}_t(t,\xi)|^2. \end{aligned} \tag{3.8}$$

*Proof.* First, to estimate $\widehat{\chi^\varepsilon}$, we take the trace of $(3.5)_3$ and use $\operatorname{tr}(\widehat{\sigma(n^\varepsilon)}) = 0$ to obtain

$$\varepsilon \operatorname{tr}(\widehat{\chi^\varepsilon})_t + \operatorname{tr}(\widehat{\chi^\varepsilon}) + \lambda d(\mathrm{i}\xi \cdot \widehat{n^\varepsilon}) = \varepsilon^{1-\alpha}\lambda d(\mathrm{i}\xi \cdot \widehat{u})_t + \operatorname{tr}(\widehat{F}_2). \tag{3.9}$$

Multiplying (3.9) by $\operatorname{tr}(\overline{\widehat{\chi^\varepsilon}})$ and taking the real part, we obtain

$$\begin{aligned} \operatorname{Re}(\mathrm{i}\xi \cdot \widehat{n^\varepsilon} \operatorname{tr}(\overline{\widehat{\chi^\varepsilon}})) =& -(\lambda d)^{-1}\left(\frac{\varepsilon}{2}\partial_t(|\operatorname{tr}(\widehat{\chi^\varepsilon})|^2) + |\operatorname{tr}(\widehat{\chi^\varepsilon})|^2\right) \\ &+ \varepsilon^{1-\alpha}\operatorname{Re}\left((\mathrm{i}\xi \cdot \widehat{u})_t \operatorname{tr}(\overline{\widehat{\chi^\varepsilon}})\right) + (\lambda d)^{-1}\operatorname{Re}\left(\operatorname{tr}(\widehat{F}_2)\operatorname{tr}(\overline{\widehat{\chi^\varepsilon}})\right). \end{aligned} \tag{3.10}$$

On the other hand, taking the Frobenius product of $(3.5)_3$ with $\frac{1}{2\nu}|\xi|^{2k}\overline{\widehat{\chi^\varepsilon}}$ and then taking the real part, we obtain

$$\begin{aligned} &\frac{\varepsilon}{4\nu}\frac{d}{dt}(|\xi|^{2k}|\widehat{\chi^\varepsilon}|^2) + \frac{1}{2\nu}|\xi|^{2k}|\widehat{\chi^\varepsilon}|^2 + |\xi|^{2k}\operatorname{Re}\left(\mathrm{i}\widehat{n^\varepsilon}\xi^T : \overline{\widehat{\chi^\varepsilon}}\right) + \frac{1}{2\nu}\left(\lambda - \frac{2\nu}{d}\right)|\xi|^{2k}\operatorname{Re}\left((\mathrm{i}\xi \cdot \widehat{n^\varepsilon})\operatorname{tr}(\overline{\widehat{\chi^\varepsilon}})\right) \\ =& -\frac{1}{2\nu}\varepsilon^{1-\alpha}|\xi|^{2k}\operatorname{Re}\left(\widehat{\pi}_t : \overline{\widehat{\chi^\varepsilon}}\right) + \frac{1}{2\nu}|\xi|^{2k}\operatorname{Re}\left(\widehat{F}_2 : \overline{\widehat{\chi^\varepsilon}}\right). \end{aligned} \tag{3.11}$$

Let $\omega_1 = \frac{1}{\lambda d^2} - \frac{1}{2\nu d}$. Multiplying (3.10) by $|\xi|^{2k}$ and substituting the result into (3.11), we obtain

$$\frac{d}{dt}\left\{\frac{\varepsilon}{4\nu}|\xi|^{2k}|\widehat{\chi^\varepsilon}|^2 + \frac{1}{2}\omega_1\varepsilon|\xi|^{2k}|\operatorname{tr}(\widehat{\chi^\varepsilon})|^2\right\} + \frac{1}{2\nu}|\xi|^{2k}|\widehat{\chi^\varepsilon}|^2 + \omega_1|\xi|^{2k}|\operatorname{tr}(\widehat{\chi^\varepsilon})|^2 + |\xi|^{2k}\operatorname{Re}\left(\mathrm{i}\widehat{n^\varepsilon}\xi^T : \overline{\widehat{\chi^\varepsilon}}\right) = F_3, \tag{3.12}$$

where

$$\begin{aligned} F_3 =& -(2\nu)^{-1}\varepsilon^{1-\alpha}|\xi|^{2k}\operatorname{Re}\left(\widehat{\pi}_t : \overline{\widehat{\chi^\varepsilon}}\right) + \varepsilon^{1-\alpha}\omega_1\lambda d|\xi|^{2k}\operatorname{Re}\left((\mathrm{i}\xi \cdot \widehat{u_t})\operatorname{tr}(\overline{\widehat{\chi^\varepsilon}})\right) \\ &+ (2\nu)^{-1}|\xi|^{2k}\operatorname{Re}\left(\widehat{F}_2 : \overline{\widehat{\chi^\varepsilon}}\right) + \left(\frac{1}{\lambda d^2} - \frac{1}{2\nu d}\right)|\xi|^{2k}\operatorname{Re}\left(\operatorname{tr}(\widehat{F}_2)\operatorname{tr}(\overline{\widehat{\chi^\varepsilon}})\right). \end{aligned}$$

Next, multiplying $(3.5)_1$ by $p'(1)|\xi|^{2k}\overline{\widehat{z^\varepsilon}}$ and $(3.5)_2$ by $|\xi|^{2k}\overline{\widehat{n^\varepsilon}}$, adding the resulting identities to (3.12) and then taking the real part, we obtain

$$\begin{aligned} &\frac{1}{2}|\xi|^{2k}\frac{d}{dt}\left\{p'(1)|\widehat{z^\varepsilon}|^2 + |\widehat{n^\varepsilon}|^2 + \frac{\varepsilon}{2\nu}|\widehat{\chi^\varepsilon}|^2 + \varepsilon\omega_1|\operatorname{tr}(\widehat{\chi^\varepsilon})|^2\right\} + \frac{1}{2\nu}|\xi|^{2k}|\widehat{\chi^\varepsilon}|^2 + \omega_1|\xi|^{2k}|\operatorname{tr}(\widehat{\chi^\varepsilon})|^2 \\ =&|\xi|^{2k}\operatorname{Re}(\overline{\widehat{n^\varepsilon}} \cdot \widehat{F}_1) + \widehat{F}_3. \end{aligned} \tag{3.13}$$

Finally, we derive the dissipation estimates for $\widehat{n^\varepsilon}$ and $\widehat{z^\varepsilon}$. To estimate the dissipation of $\widehat{n^\varepsilon}$, we take the Frobenius product of $(3.5)_3$ with $\mathrm{i}|\xi|^{2k-2}\overline{\widehat{n^\varepsilon}}\xi^T$ and then take the real part. This gives

$$
\begin{aligned}
&-\{\nu\widehat{\sigma(n^\varepsilon)}+\lambda(\mathrm{i}\xi\cdot\widehat{n^\varepsilon})I_d\}:\mathrm{i}|\xi|^{2k-2}\overline{\widehat{n^\varepsilon}}\xi^T=\nu|\xi|^{2k}|\widehat{n^\varepsilon}|^2+\mu|\xi|^{2k-2}|\xi\cdot\widehat{n^\varepsilon}|^2\\
=&\varepsilon\operatorname{Re}(\widehat{\chi^\varepsilon_t}:\mathrm{i}|\xi|^{2k-2}\overline{\widehat{n^\varepsilon}}\xi^T)+\operatorname{Re}(\widehat{\chi^\varepsilon}:\mathrm{i}|\xi|^{2k-2}\overline{\widehat{n^\varepsilon}}\xi^T)+\varepsilon^{1-\alpha}\operatorname{Re}(\widehat{\pi}_t:\mathrm{i}|\xi|^{2k-2}\overline{\widehat{n^\varepsilon}}\xi^T)-\operatorname{Re}(\widehat{F}_2:\mathrm{i}|\xi|^{2k-2}\overline{\widehat{n^\varepsilon}}\xi^T).
\end{aligned}
\tag{3.14}
$$

For the first term on the right-hand side of (3.14), $(3.5)_2$ yields

$$
\begin{aligned}
\varepsilon\operatorname{Re}(\widehat{\chi^\varepsilon_t}:\mathrm{i}|\xi|^{2k-2}\overline{\widehat{n^\varepsilon}}\xi^T)=&\,\varepsilon|\xi|^{2k-2}\operatorname{Im}(\overline{\widehat{\chi^\varepsilon}}:\widehat{n^\varepsilon}\xi^T)_t\\
&+\varepsilon|\xi|^{2k-2}\operatorname{Im}(\mathrm{i}\widehat{\chi^\varepsilon}\xi\xi^T:\overline{\widehat{\chi^\varepsilon}})+p'(1)\varepsilon|\xi|^{2k-2}\operatorname{Im}(\mathrm{i}\widehat{z^\varepsilon}\xi\xi^T:\overline{\widehat{\chi^\varepsilon}})-\varepsilon|\xi|^{2k-2}\xi^T\operatorname{Im}(\overline{\widehat{\chi^\varepsilon}}\widehat{F}_1).
\end{aligned}
\tag{3.15}
$$

For the second term on the right-hand side of (3.14), we multiply $(3.5)_2$ by $|\xi|^{2k-2}\overline{\widehat{n^\varepsilon}}$, take the real part and use $\operatorname{Re}(\mathrm{i}\widehat{z^\varepsilon}\xi\cdot\overline{\widehat{n^\varepsilon}})=-\operatorname{Re}(\widehat{z^\varepsilon}\overline{\mathrm{i}\xi\cdot\widehat{n^\varepsilon}})$ to obtain

$$
\begin{aligned}
&\operatorname{Re}(\widehat{\chi^\varepsilon}:\mathrm{i}|\xi|^{2k-2}\overline{\widehat{n^\varepsilon}}\xi^T)=|\xi|^{2k-2}\operatorname{Re}(\mathrm{i}\widehat{\chi^\varepsilon}\xi\cdot\overline{\widehat{n^\varepsilon}})\\
=&-\frac{1}{2}|\xi|^{2k-2}|\widehat{n^\varepsilon}|^2_t+p'(1)|\xi|^{2k-2}\operatorname{Re}(\widehat{z^\varepsilon}\overline{\mathrm{i}\xi\cdot\widehat{n^\varepsilon}})+|\xi|^{2k-2}\operatorname{Re}(\widehat{F}_1\cdot\overline{\widehat{n^\varepsilon}}).
\end{aligned}
\tag{3.16}
$$

Taking the complex conjugate of $(3.5)_1$, multiplying it by $|\xi|^{2k-2}\widehat{z^\varepsilon}$ and taking the real part, we obtain

$$
|\xi|^{2k-2}\operatorname{Re}(\widehat{z^\varepsilon}\overline{\mathrm{i}\xi\cdot\widehat{n^\varepsilon}})=-\frac{1}{2}|\xi|^{2k-2}(|\widehat{z^\varepsilon}|^2)_t.
\tag{3.17}
$$

Combining (3.16) and (3.17), we have

$$
\operatorname{Re}(\widehat{\chi^\varepsilon}:\mathrm{i}|\xi|^{2k-2}\overline{\widehat{n^\varepsilon}}\xi^T)=-\frac{1}{2}|\xi|^{2k-2}(|\widehat{n^\varepsilon}|^2+p'(1)|\widehat{z^\varepsilon}|^2)_t+|\xi|^{2k-2}\operatorname{Re}(\widehat{F}_1\cdot\overline{\widehat{n^\varepsilon}}).
\tag{3.18}
$$

Substituting (3.15) and (3.18) into (3.14), we obtain

$$
\begin{aligned}
&\nu|\xi|^{2k}|\widehat{n^\varepsilon}|^2+\mu|\xi|^{2k-2}|\xi\cdot\widehat{n^\varepsilon}|^2\\
=&-\frac{1}{2}|\xi|^{2k-2}(|\widehat{n^\varepsilon}|^2+p'(1)|\widehat{z^\varepsilon}|^2)_t+\varepsilon|\xi|^{2k-2}\operatorname{Im}(\overline{\widehat{\chi^\varepsilon}}:\widehat{n^\varepsilon}\xi^T)_t+\varepsilon|\xi|^{2k-2}\operatorname{Im}(\mathrm{i}\widehat{\chi^\varepsilon}\xi\xi^T:\overline{\widehat{\chi^\varepsilon}})\\
&+p'(1)\varepsilon|\xi|^{2k-2}\operatorname{Im}(\mathrm{i}\widehat{z^\varepsilon}\xi\xi^T:\overline{\widehat{\chi^\varepsilon}})-\varepsilon|\xi|^{2k-2}\xi^T\operatorname{Im}(\overline{\widehat{\chi^\varepsilon}}\widehat{F}_1)\\
&+|\xi|^{2k-2}\operatorname{Re}(\widehat{F}_1\cdot\overline{\widehat{n^\varepsilon}})+\varepsilon^{1-\alpha}|\xi|^{2k-2}\operatorname{Re}(\widehat{\pi}_t:\mathrm{i}\overline{\widehat{n^\varepsilon}}\xi^T)-|\xi|^{2k-2}\operatorname{Re}(\widehat{F}_2:\mathrm{i}\overline{\widehat{n^\varepsilon}}\xi^T).
\end{aligned}
\tag{3.19}
$$

To estimate the dissipation of $\widehat{z^\varepsilon}$, we multiply $(3.5)_2$ by $\mathrm{i}|\xi|^{2k-2}\xi\overline{\widehat{z^\varepsilon}}$, take the real part and use the identity

$$
\operatorname{Re}(\mathrm{i}\widehat{n^\varepsilon_t}\overline{\widehat{z^\varepsilon}})=\operatorname{Im}(\overline{\widehat{n^\varepsilon}}\widehat{z^\varepsilon})_t-\operatorname{Im}(\overline{\widehat{n^\varepsilon}}\widehat{z^\varepsilon_t}).
$$

We then obtain

$$
\begin{aligned}
p'(1)|\xi|^{2k}|\widehat{z^\varepsilon}|^2=&\operatorname{Re}(\widehat{n^\varepsilon_t}\cdot\mathrm{i}|\xi|^{2k-2}\xi\overline{\widehat{z^\varepsilon}})+\operatorname{Re}(\mathrm{i}\widehat{\chi^\varepsilon}\xi\cdot\mathrm{i}|\xi|^{2k-2}\xi\overline{\widehat{z^\varepsilon}})-\operatorname{Re}(\widehat{F}_1\cdot\mathrm{i}|\xi|^{2k-2}\xi\overline{\widehat{z^\varepsilon}})\\
=&|\xi|^{2k-2}\xi\cdot\operatorname{Im}(\overline{\widehat{n^\varepsilon}}\widehat{z^\varepsilon})_t+|\xi|^{2k-2}|\xi\cdot\widehat{n^\varepsilon}|^2-|\xi|^{2k-2}\operatorname{Re}(\widehat{\chi^\varepsilon}\xi\cdot\xi\overline{\widehat{z^\varepsilon}})-|\xi|^{2k-2}\xi\cdot\operatorname{Re}(\mathrm{i}\widehat{F}_1\overline{\widehat{z^\varepsilon}}).
\end{aligned}
\tag{3.20}
$$

Choose $\epsilon_4$ to be a sufficiently small positive constant. Multiplying (3.20) by $\epsilon_4$ and adding the resulting identity to (3.13) and (3.19), we obtain

$$\frac{\partial(G(t,\xi)^2)}{\partial t} + D(t,\xi) = \widetilde{Q}, \tag{3.21}$$

where

$$\begin{aligned}
G(t,\xi)^2 =& \frac{1}{2}|\xi|^{2k}\left\{p'(1)|\widehat{z^\varepsilon}|^2 + |\widehat{n^\varepsilon}|^2 + \frac{\varepsilon}{2\nu}|\widehat{\chi^\varepsilon}|^2 + \varepsilon\omega_1|\operatorname{tr}(\widehat{\chi^\varepsilon})|^2\right\}\\
&+ \frac{1}{2}|\xi|^{2k-2}(p'(1)|\widehat{z^\varepsilon}|^2 + |\widehat{n^\varepsilon}|^2) - \varepsilon|\xi|^{2k-2}\operatorname{Im}(\overline{\widehat{\chi^\varepsilon}} : \widehat{n^\varepsilon}\xi^T) - \epsilon_4|\xi|^{2k-2}\xi\cdot\operatorname{Im}(\overline{\widehat{n^\varepsilon}}\widehat{z^\varepsilon}),\\
D(t,\xi) =& \frac{1}{2\nu}|\xi|^{2k}|\widehat{\chi^\varepsilon}|^2 - \varepsilon|\xi|^{2k-2}\operatorname{Im}(\mathrm{i}\widehat{\chi^\varepsilon}\xi\xi^T : \overline{\widehat{\chi^\varepsilon}}) + \omega_1|\xi|^{2k}|\operatorname{tr}(\widehat{\chi^\varepsilon})|^2 + \nu|\xi|^{2k}|\widehat{n^\varepsilon}|^2 - \epsilon_4|\xi|^{2k-2}|\xi\cdot\widehat{n^\varepsilon}|^2\\
&+ p'(1)\epsilon_4|\xi|^{2k}|\widehat{z^\varepsilon}|^2 + \mu|\xi|^{2k-2}|\xi\cdot\widehat{n^\varepsilon}|^2 - p'(1)\varepsilon|\xi|^{2k-2}\operatorname{Im}(\mathrm{i}\widehat{z^\varepsilon}\xi\xi^T : \overline{\widehat{\chi^\varepsilon}})\\
&+ \epsilon_4|\xi|^{2k-2}\operatorname{Re}(\widehat{\chi^\varepsilon}\xi\cdot\xi\overline{\widehat{z^\varepsilon}}),\\
\sim& |\xi|^{2k}\left(|\widehat{z^\varepsilon}(t,\xi)|^2 + |\widehat{n^\varepsilon}(t,\xi)|^2 + |\widehat{\chi^\varepsilon}(t,\xi)|^2\right),\\
\widetilde{Q}(t,\xi) =& \varepsilon^{1-\alpha}|\xi|^{2k-2}\operatorname{Re}(\widehat{\pi}_t : \mathrm{i}\overline{\widehat{n^\varepsilon}}\xi^T) + F_3 - |\xi|^{2k-2}\operatorname{Re}(\widehat{F}_2 : \mathrm{i}\overline{\widehat{n^\varepsilon}}\xi^T) + |\xi|^{2k-2}\operatorname{Re}(\widehat{F}_1\cdot\overline{\widehat{n^\varepsilon}})\\
&+ |\xi|^{2k}\operatorname{Re}(\widehat{F}_1\cdot\overline{\widehat{n^\varepsilon}}) - \varepsilon|\xi|^{2k-2}\xi^T\cdot\operatorname{Im}(\overline{\widehat{\chi^\varepsilon}}\widehat{F}_1) - \epsilon_4|\xi|^{2k-2}\xi\cdot\operatorname{Re}(\mathrm{i}\widehat{F}_1\overline{\widehat{z^\varepsilon}}).
\end{aligned}$$

Since $|\operatorname{tr}(\widehat{\chi^\varepsilon})|^2 \le d|\widehat{\chi^\varepsilon}|^2$, by choosing $\epsilon_4$ and $\varepsilon$ sufficiently small, it is readily seen that $G(t,\xi)$ satisfies (3.6). We now estimate the source terms in (3.21). By (3.6) and $0 < \varepsilon \le 1$,

$$\begin{aligned}
&\left(|\xi|^{2k-2} + |\xi|^{2k}\right)|\widehat{F}_1|\,|\widehat{n^\varepsilon}| + |\xi|^{2k-1}|\widehat{F}_1|\,|\widehat{z^\varepsilon}| + \varepsilon|\xi|^{2k-1}|\widehat{F}_1|\,|\widehat{\chi^\varepsilon}| + \varepsilon^{1-\alpha}|\xi|^{2k-1}|\widehat{\pi}_t|\,|\widehat{n^\varepsilon}|\\
&\lesssim \left(|\xi|^{k-1} + |\xi|^k\right)|\widehat{F}_1|\,G + \varepsilon^{1-\alpha}|\xi|^k|\widehat{\pi}_t|\,G \lesssim A_k G.
\end{aligned} \tag{3.22}$$

For the remaining terms, we apply Young's inequality to obtain

$$\begin{aligned}
F_3 + |\xi|^{2k-1}|\widehat{F}_2|\,|\widehat{n^\varepsilon}| &\lesssim \epsilon D + \left(|\xi|^{2k-2} + |\xi|^{2k}\right)|\widehat{F}_2|^2 + \varepsilon^{2-2\alpha}|\xi|^{2k}|\widehat{\pi}_t|^2\\
&\lesssim \epsilon D + B_k^2.
\end{aligned} \tag{3.23}$$

Combining (3.22) and (3.23), we obtain

$$\widetilde{Q} \lesssim \epsilon D + A_k G + B_k^2. \tag{3.24}$$

Observing that $D \gtrsim \lambda_\xi G^2$ and substituting (3.24) into (3.21), we obtain (3.7), with suitable constants $c, C_0 > 0$. □

With the Fourier energy inequality (3.7) established, we now prove Proposition 3.1.

*Proof of Proposition 3.1.* Fix $T > 0$. We first establish (3.1). Set $G_0(\xi) = G(0,\xi)$ and define the auxiliary functions $Z(t,\xi)$ and $W(t,\xi)$ by

$$\begin{aligned}
Z(t,\xi) &= e^{-c\lambda_\xi t}G_0(\xi) + \frac{C_0}{2}\int_0^t e^{-c\lambda_\xi(t-\tau)}A_k(\tau,\xi)\,d\tau,\\
W(t,\xi)^2 &= C_0\int_0^t e^{-2c\lambda_\xi(t-\tau)}B_k(\tau,\xi)^2\,d\tau.
\end{aligned}$$

Applying Lemma A.2 in the Appendix to the differential inequality (3.7), we obtain

$$G(t,\xi)\le Z(t,\xi)+W(t,\xi),\qquad 0\le t\le T.$$

Taking the $L^2_\xi$-norm and using (3.6) together with Plancherel's theorem, we obtain

$$\|G(t)\|_{L^2_\xi}\sim\|D^{k-1}(z^\varepsilon(t),n^\varepsilon(t))\|_1+\sqrt{\varepsilon}\|D^k\chi^\varepsilon(t)\|.$$

By Lemma A.2, Minkowski's inequality and (3.6), we obtain

$$\|D^{k-1}(z^\varepsilon,n^\varepsilon)(t)\|_{H^1}+\sqrt{\varepsilon}\,\|D^k\chi^\varepsilon(t)\|\lesssim\|G(t)\|_{L^2_\xi}\lesssim\sum_{j=0}^{6}Q_j(t), \tag{3.25}$$

where

$$\begin{aligned}
Q_0(t)&:=\left\|e^{-c\lambda_\xi t}G(0,\xi)\right\|_{L^2_\xi},\\
Q_1(t)&:=\varepsilon^{1-\alpha}\int_0^t\left\|e^{-c\lambda_\xi(t-\tau)}|\xi|^k\widehat{\pi}_t(\tau)\right\|_{L^2_\xi}\,d\tau,\\
Q_2(t)&:=\varepsilon^{1-\alpha}\left(\int_0^t\int_{\mathbb{R}^d}e^{-2c\lambda_\xi(t-\tau)}|\xi|^{2k}\left|\widehat{\pi}_t(\tau,\xi)\right|^2d\xi\,d\tau\right)^{1/2},\\
Q_3(t)&:=\int_0^t\left\|e^{-c\lambda_\xi(t-\tau)}|\xi|^{k-1}\widehat{F_1}(\tau)\right\|_{L^2_\xi}\,d\tau,\\
Q_4(t)&:=\int_0^t\left\|e^{-c\lambda_\xi(t-\tau)}|\xi|^{k}\widehat{F_1}(\tau)\right\|_{L^2_\xi}\,d\tau,\\
Q_5(t)&:=\left(\int_0^t\int_{\mathbb{R}^d}e^{-2c\lambda_\xi(t-\tau)}|\xi|^{2k-2}\left|\widehat{F_2}(\tau,\xi)\right|^2d\xi\,d\tau\right)^{1/2},\\
Q_6(t)&:=\left(\int_0^t\int_{\mathbb{R}^d}e^{-2c\lambda_\xi(t-\tau)}|\xi|^{2k}\left|\widehat{F_2}(\tau,\xi)\right|^2d\xi\,d\tau\right)^{1/2}.
\end{aligned}$$

For $k=1,\cdots,[\frac{d}{2}]+2$, in view of the definition of $Y(t)$, it suffices to establish the following estimates:

$$\begin{aligned}
&Q_0\lesssim(1+t)^{-\frac{d+2k}{4}},\\
&Q_j\lesssim\varepsilon^{1-\alpha}(1+t)^{-\frac{d+2k}{4}},\quad j=1,2,\\
&Q_j\lesssim(1+t)^{-\frac{d+2k}{4}}(C(\epsilon)+\epsilon Y(t)+\varepsilon^\alpha Y(t)^2),\quad j=3,\cdots,6.
\end{aligned}\tag{3.26}$$

Among the above estimates, the cases $j=0,1,4$ and 6 are the most representative. We prove only these four cases below.

We first estimate $Q_0$. For all $k=1,\dots,[\frac{d}{2}]+2$, we have

$$Q_0\lesssim\left(\int_{|\xi|\ge1}e^{-C\lambda_\xi t}\left(|\xi|^{2k-2}(1+|\xi|^2)(|\widehat{z_0^\varepsilon}|^2+|\widehat{n_0^\varepsilon}|^2)+\varepsilon|\xi|^{2k}|\widehat{\chi_0^\varepsilon}|^2\right)d\xi\right)^{\frac12}$$

$$
\begin{aligned}
&+\left(\int_{|\xi|\le 1} e^{-C\lambda_\xi t}\left(|\xi|^{2k-2}(1+|\xi|^2)(|\widehat{z_0^\varepsilon}|^2+|\widehat{n_0^\varepsilon}|^2)+\varepsilon|\xi|^{2k}|\widehat{\chi_0^\varepsilon}|^2\right)d\xi\right)^{\frac12}\\
&=: Q_0^1+Q_0^2.
\end{aligned}
$$

From the assumed condition (1.13), we derive the fact that $(\widehat{z_0^\varepsilon}(0),\widehat{n_0^\varepsilon}(0))=(0,0)$. Combining this with

$$
\|\nabla_\xi \hat f(t)\|_{L^\infty_\xi}\le \int_{\mathbb{R}^d}(1+|x|)|f(t,x)|\,dx,
$$

we obtain

$$
\begin{aligned}
Q_0^1&\lesssim\left(\int_{|\xi|\ge 1} e^{-Ct}\left(|\xi|^{2k-2}(1+|\xi|^2)(|\widehat{n_0^\varepsilon}|^2+|\widehat{z_0^\varepsilon}|^2)+\varepsilon|\xi|^{2k}|\widehat{\chi_0^\varepsilon}|^2\right)d\xi\right)^{\frac12}\\
&\lesssim e^{-Ct},
\end{aligned}
$$

and

$$
\begin{aligned}
Q_0^2&\lesssim\left(\int_{|\xi|\le 1} e^{-C|\xi|^2(1+t)}|\xi|^{2k}\left(\left(1+\frac{1}{|\xi|^2}\right)(|\widehat{n_0^\varepsilon}|^2+|\widehat{z_0^\varepsilon}|^2)+\varepsilon|\widehat{\chi_0^\varepsilon}|^2\right)d\xi\right)^{\frac12}\\
&\lesssim\left(\int_{|\xi|\le 1} e^{-C|\xi|^2(1+t)}|\xi|^{2k}d\xi\right)^{\frac12}\\
&\lesssim (1+t)^{-\frac{d+2k}{4}}.
\end{aligned}
$$

Consequently, (3.26) for $j=0$ holds.

Next, we estimate $Q_1$. Using $(1.3)_2$, we observe that

$$
\begin{aligned}
u_t=&-h'(1)\nabla(\rho-1)+\nu\Delta u+\mu\nabla(\nabla\cdot u)\\
&-u\cdot\nabla u-\nabla\left(h(\rho)-h'(1)(\rho-1)\right)-\rho^{-1}\left(\nu\Delta u+\mu\nabla(\nabla\cdot u)\right)(\rho-1).
\end{aligned}
$$

For $|\xi|\le 1$, we have

$$
\begin{aligned}
|\widehat{\pi}_t|\lesssim&\underbrace{|\xi|^2\left|(\widehat{\rho-1},\widehat{u})\right|}_{linear\quad terms}\\
&+\underbrace{|\xi|\left(|\mathcal{F}(u\cdot\nabla u)|+|\xi||\mathcal{F}\left(h(\rho)-h'(1)(\rho-1)\right)|+|\mathcal{F}\left(\rho^{-1}\left(\nu\Delta u+\mu\nabla(\nabla\cdot u)\right)(\rho-1)\right)|\right)}_{nonlinear\quad terms}\\
=:&|\xi|^2\left|(\widehat{\rho-1},\widehat{u})\right|+|\xi|\mathcal{J}(t,\xi).
\end{aligned}
$$

We apply Lemma 3.2 to the terms involving the above linear part and use the Fourier transform property $\|\hat f(t)\|_{L^\infty}\lesssim\|f(t)\|_{L^1}$ for the remaining terms involving the nonlinear part. It follows

that

$$
\begin{aligned}
Q_1(t) &\lesssim \varepsilon^{1-\alpha} \int_0^t \left( \int_{|\xi|\geq 1} e^{-C(t-\tau)} |\xi|^{2k} |\widehat{\pi}_t(\tau)|^2 \, d\xi \right)^{1/2} d\tau \\
&\quad + \varepsilon^{1-\alpha} \int_{\frac{t}{2}}^t \left( \int_{|\xi|\leq 1} e^{-C|\xi|^2(1+t-\tau)} |\xi|^{2k} |\widehat{\pi}_t(\tau)|^2 \, d\xi \right)^{1/2} d\tau \\
&\quad + \varepsilon^{1-\alpha} \int_0^{\frac{t}{2}} \left( \int_{|\xi|\leq 1} e^{-C|\xi|^2(1+t-\tau)} |\xi|^{2k+2} |\mathcal{J}(\tau)|^2 \, d\xi \right)^{1/2} d\tau \\
&\quad + \varepsilon^{1-\alpha} \int_0^{\frac{t}{2}} \left( \int_{|\xi|\leq 1} e^{-C|\xi|^2(1+t-\tau)} |\xi|^{2k+4} |(\widehat{\rho - 1}, \widehat{u})(\tau)|^2 \, d\xi \right)^{1/2} d\tau \\
&\lesssim \varepsilon^{1-\alpha} \left( \int_0^{\frac{t}{2}} + \int_{\frac{t}{2}}^t \right) e^{-C(t-\tau)} \|D^{k+1} u_t(\tau)\| \, d\tau + \varepsilon^{1-\alpha} \int_{\frac{t}{2}}^t \|D^{k+1} u_t(\tau)\| \, d\tau \\
&\quad + \varepsilon^{1-\alpha} \int_0^{\frac{t}{2}} \Big( \|u(\tau)\| \|Du(\tau)\| + \|\rho(\tau) - 1\| \|D\rho(\tau)\| \\
&\qquad\qquad + \|\rho(\tau) - 1\| \|D^2 u(\tau)\| \Big) \left( \int_{|\xi|\leq 1} e^{-C|\xi|^2(1+t-\tau)} |\xi|^{2k+2} \, d\xi \right)^{1/2} d\tau \\
&\quad + \varepsilon^{1-\alpha} \sup_{\tau\geq 0} \|(\widehat{\rho - 1}, \widehat{u})(\tau)\|_{L^\infty_\xi} \int_0^{\frac{t}{2}} \left( \int_{|\xi|\leq 1} e^{-C|\xi|^2(1+t-\tau)} |\xi|^{2k+4} \, d\xi \right)^{1/2} d\tau \\
&\lesssim \varepsilon^{1-\alpha} (1+t)^{-\frac{d+2k}{4}}.
\end{aligned} \tag{3.27}
$$

This proves the second estimate in (3.26) for $j = 1$.

We now estimate $Q_4$. By Minkowski's inequality and the conservative structure of $F_1$, we obtain

$$
\begin{aligned}
Q_4(t) &\lesssim \int_0^t \left\| e^{-\frac{c|\xi|^2(t-\tau)}{1+|\xi|^2}} |\xi|^{k+1} \mathcal{F}\left( \frac{m \otimes n^\varepsilon}{\rho^\varepsilon} \right) \right\|_{L^2_\xi} d\tau \\
&\quad + \int_0^t \left\| e^{-\frac{c|\xi|^2(t-\tau)}{1+|\xi|^2}} |\xi|^{k+1} \mathcal{F}\left( \frac{n^\varepsilon \otimes m}{\rho^\varepsilon} \right) \right\|_{L^2_\xi} d\tau \\
&\quad + \varepsilon^\alpha \int_0^t \left\| e^{-\frac{c|\xi|^2(t-\tau)}{1+|\xi|^2}} |\xi|^{k+1} \mathcal{F}\left( \frac{n^\varepsilon \otimes n^\varepsilon}{\rho^\varepsilon} \right) \right\|_{L^2_\xi} d\tau \\
&\quad + \int_0^t \left\| e^{-\frac{c|\xi|^2(t-\tau)}{1+|\xi|^2}} |\xi|^{k+1} \mathcal{F}\left( \frac{z^\varepsilon m \otimes m}{\rho \rho^\varepsilon} \right) \right\|_{L^2_\xi} d\tau \\
&\quad + \int_0^t \left\| e^{-\frac{c|\xi|^2(t-\tau)}{1+|\xi|^2}} |\xi|^{k+1} \mathcal{F}\left( \big(p'(\rho, \rho^\varepsilon) - p'(1)\big) z^\varepsilon \right) \right\|_{L^2_\xi} d\tau \\
&=: \sum_{j=1}^5 K_j.
\end{aligned}
$$

We first estimate $K_1$. Throughout the following estimates, $1 \leq k \leq [\frac{d}{2}] + 2$ and $0 \leq \tau \leq t$.

The uniform energy estimate and the definition of $Y(t)$ imply

$$
\begin{aligned}
&\|(z^\varepsilon(t), n^\varepsilon(t))\|_{L^\infty} \lesssim Y(t)(1+t)^{-\frac{d+1}{2}},\\
&\|D^k(z^\varepsilon(t), n^\varepsilon(t))\| \lesssim \begin{cases} Y(t)(1+t)^{-\frac{d+2k+2}{4}}, & k=0,\cdots,[\frac{d}{2}]+1,\\ Y(t)(1+t)^{-\frac{d+2[\frac{d}{2}]+4}{4}}, & k=[\frac{d}{2}]+2.\end{cases}
\end{aligned}
$$

A frequency decomposition gives

$$
\begin{aligned}
K_1 \lesssim& \int_0^t \left(\int_{|\xi|\geq 1} e^{-\frac{2c|\xi|^2(t-\tau)}{1+|\xi|^2}}|\xi|^{2k+2}\left|\mathcal{F}\left(\frac{m\otimes n^\varepsilon}{\rho^\varepsilon}\right)\right|^2 d\xi\right)^{1/2} d\tau\\
&+\int_0^t \left(\int_{|\xi|\leq 1} e^{-\frac{2c|\xi|^2(t-\tau)}{1+|\xi|^2}}|\xi|^{2k+2}\left|\mathcal{F}\left(\frac{m\otimes n^\varepsilon}{\rho^\varepsilon}\right)\right|^2 d\xi\right)^{1/2} d\tau\\
=:& K_1^1+K_1^2.
\end{aligned}
\tag{3.28}
$$

For the high-frequency contribution, the Moser-type inequalities yield

$$
\begin{aligned}
K_1^1 \lesssim& \int_0^t e^{-C(t-\tau)}\left\|D^{k+1}\left(\frac{m\otimes n^\varepsilon}{\rho^\varepsilon}\right)\right\| d\tau\\
\lesssim& \int_0^t e^{-C(t-\tau)}\Big(\|D^{k+1}n^\varepsilon\|\|m\|_{L^\infty}+\|D^{k+1}m\|\|n^\varepsilon\|_{L^\infty}\\
&\qquad +\|D^{k+1}(\rho,\varepsilon^\alpha z^\varepsilon)\|\|m\|_{L^\infty}\|n^\varepsilon\|_{L^\infty}\Big) d\tau.
\end{aligned}
$$

Noticing the Gagliardo-Nirenberg inequality

$$
\|D^{k+1}n^\varepsilon\| \lesssim \|D^k n^\varepsilon\|^{\gamma_1}\|D^s n^\varepsilon\|^{1-\gamma_1}, \quad \gamma_1=\frac{s-k-1}{s-k}
\tag{3.29}
$$

and using the uniform energy bound in Theorem 1.1 together with assumption ($H_2$), under the assumption $s\geq[\frac{d}{2}]+4$ in Theorem 1.2, when $k=1,\cdots,[\frac{d}{2}]+1$, we obtain

$$
\begin{aligned}
K_1^1 \lesssim& \int_0^t e^{-C(t-\tau)}\Big((1+\tau)^{-\frac{d+2k+2}{4}}+Y(t)^{\gamma_1}(1+\tau)^{-\frac{d}{2}-\gamma_1\frac{d+2k+2}{4}}\Big) d\tau\\
\lesssim& (1+t)^{-\frac{d+2k}{4}}\big(1+Y(t)^{\gamma_1}\big).
\end{aligned}
\tag{3.30}
$$

Here we used

$$
\frac{d}{2}+\gamma_1\frac{d+2k+2}{4}\geq\frac{d+2k}{4}.
$$

When $k=[\frac{d}{2}]+2$,

$$
\begin{aligned}
K_1^1 \lesssim& \int_0^t e^{-C(t-\tau)}\Big((1+\tau)^{-\frac{d+2k+2}{4}}+Y(t)^{\gamma_1}(1+\tau)^{-\frac{d}{2}-\gamma_1\frac{d+2k}{4}}\Big) d\tau\\
\lesssim& (1+t)^{-\frac{d+2k}{4}}\big(1+Y(t)^{\gamma_1}\big).
\end{aligned}
\tag{3.31}
$$

Here we used

$$
\frac{d}{2}+\gamma_1\frac{d+2[\frac{d}{2}]+4}{4}\geq\frac{d+2[\frac{d}{2}]+4}{4}.
$$

We next estimate the low-frequency contribution. Choose $\gamma \in [\frac{3}{4}, 1)$. For $0 \le \tau \le \frac{t}{2}$, since $|\xi|^k \le |\xi|^{k-1}$ for $|\xi| \le 1$, the Fourier transform bound $\|\hat{f}(t)\|_{L^\infty_\xi} \lesssim \|f(t)\|_{L^1}$ and Hölder's inequality give

$$
\begin{aligned}
&\int_0^{\frac{t}{2}} \left( \int_{|\xi|\le 1} e^{-C|\xi|^2(t-\tau)} |\xi|^{2k} \left| \mathcal{F}\left( \frac{m \otimes n^\varepsilon}{\rho^\varepsilon} \right) \right|^2 d\xi \right)^{1/2} d\tau \\
&\quad \lesssim \int_0^{\frac{t}{2}} \left\| \frac{m \otimes n^\varepsilon}{\rho^\varepsilon} \right\|_{L^1} \left( \int_{|\xi|\le 1} e^{-C|\xi|^2(1+t-\tau)} |\xi|^{2k} \, d\xi \right)^{1/2} d\tau \\
&\quad \lesssim Y(t)^\gamma \int_0^{\frac{t}{2}} (1+\tau)^{-\gamma\frac{d+2}{4} - \frac{d}{4}} (1+t-\tau)^{-\frac{d+2k}{4}} \, d\tau \\
&\quad \lesssim (1+t)^{-\frac{d+2k}{4}} Y(t)^\gamma.
\end{aligned} \tag{3.32}
$$

The last inequality in (3.32) follows from $\gamma\frac{d+2}{4} + \frac{d}{4} > 1$ for $d \ge 2$.

For $\frac{t}{2} \le \tau \le t$, we retain one power of frequency in the Fourier multiplier. For $1 \le k \le [d/2] + 2$, the Moser-type inequalities give

$$
\begin{aligned}
&\left\| D^{k-1} \left( \frac{m \otimes n^\varepsilon}{\rho^\varepsilon} \right) \right\| \\
\lesssim& \|D^{k-1} n^\varepsilon\| \|m\|_{L^\infty} + \|D^{k-1} m\| \|n^\varepsilon\|_{L^\infty} + \|D^{k-1}(\rho - 1, \varepsilon^\alpha z^\varepsilon)\| \|m\|_{L^\infty} \|n^\varepsilon\|_{L^\infty} \\
\lesssim& Y(t)^\gamma \left( (1+\tau)^{-\frac{d}{2} - \gamma\frac{d+2k}{4}} + (1+\tau)^{-\frac{d+2k-2}{4} - \gamma\frac{d+1}{2}} \right) \\
\lesssim& Y(t)^\gamma (1+\tau)^{-\frac{d+2k+2}{4}}.
\end{aligned} \tag{3.33}
$$

This choice of $\gamma$ also ensures that

$$
\frac{d}{2} + \gamma\frac{d+2k}{4} \ge \frac{d+2k+2}{4}, \qquad \frac{d+2k-2}{4} + \gamma\frac{d+1}{2} \ge \frac{d+2k+2}{4}.
$$

Using (3.33), we obtain

$$
\begin{aligned}
&\int_{\frac{t}{2}}^{t} \left( \int_{|\xi|\le 1} e^{-C|\xi|^2(t-\tau)} |\xi|^{2k} \left| \mathcal{F}\left( \frac{m \otimes n^\varepsilon}{\rho^\varepsilon} \right) \right|^2 d\xi \right)^{1/2} d\tau \\
&\quad \lesssim \int_{\frac{t}{2}}^{t} \sup_{|\xi|\le 1} \left( e^{-C|\xi|^2(1+t-\tau)} |\xi| \right) \left\| D^{k-1} \left( \frac{m \otimes n^\varepsilon}{\rho^\varepsilon} \right) \right\| d\tau \\
&\quad \lesssim Y(t)^\gamma \int_{\frac{t}{2}}^{t} (1+t-\tau)^{-\frac{1}{2}} (1+\tau)^{-\frac{d+2k+2}{4}} \, d\tau \\
&\quad \lesssim (1+t)^{-\frac{d+2k}{4}} Y(t)^\gamma.
\end{aligned} \tag{3.34}
$$

Combining (3.32) and (3.34), we obtain

$$
K_1^2 \lesssim (1+t)^{-\frac{d+2k}{4}} Y(t)^\gamma.
$$

Therefore, (3.28), (3.30), (3.31) and Young's inequality imply that

$$
K_1 \lesssim (1+t)^{-\frac{d+2k}{4}} \left( 1 + Y(t)^\gamma + Y(t)^{\gamma_1} \right)
$$

$$\lesssim (1+t)^{-\frac{d+2k}{4}}\big(C(\epsilon)+\epsilon Y(t)\big).$$

The same argument applies to $K_2$ and $K_4$. Thus, after adjusting $\epsilon$, we have

$$K_1+K_2+K_4 \leq (1+t)^{-\frac{d+2k}{4}}\big(C(\epsilon)+\epsilon Y(t)\big). \tag{3.35}$$

For the quadratic term $K_3$, we retain its factor $\varepsilon^\alpha$. The frequency and time decompositions give

$$\begin{aligned}K_3 &\lesssim \varepsilon^\alpha \int_0^t e^{-C(t-\tau)} \left\| D^{k+1}\left(\frac{n^\varepsilon\otimes n^\varepsilon}{\rho^\varepsilon}\right)\right\| d\tau \\ &\quad + \varepsilon^\alpha \int_0^{\frac{t}{2}} (1+t-\tau)^{-\frac{d+2k+2}{4}} \left\|\frac{n^\varepsilon\otimes n^\varepsilon}{\rho^\varepsilon}\right\|_{L^1} d\tau \\ &\quad + \varepsilon^\alpha \int_{\frac{t}{2}}^t (1+t-\tau)^{-\frac{1}{2}} \left\| D^{k-1}\left(\frac{n^\varepsilon\otimes n^\varepsilon}{\rho^\varepsilon}\right)\right\| d\tau.\end{aligned}$$

Using the assumptions $(H_1)$ and $(H_2)$, (3.29), the Moser-type inequalities and the definition of $Y(t)$, we have

$$\begin{aligned}\left\|\frac{n^\varepsilon\otimes n^\varepsilon}{\rho^\varepsilon}\right\|_{L^1} &\lesssim Y(t)^2(1+\tau)^{-\frac{d+2}{2}}, \\ \left\| D^{k-1}\left(\frac{n^\varepsilon\otimes n^\varepsilon}{\rho^\varepsilon}\right)\right\| &\lesssim Y(t)^2(1+\tau)^{-\frac{d+1}{2}-\frac{d+2k}{4}}, \\ \left\| D^{k+1}\left(\frac{n^\varepsilon\otimes n^\varepsilon}{\rho^\varepsilon}\right)\right\| &\lesssim \big(Y(t)^{1+\gamma_1}+Y(t)^2\big)(1+\tau)^{-\frac{d+2k}{4}}.\end{aligned}$$

Combining this with the estimate for $K_3$, we arrive at

$$\begin{aligned}K_3 &\lesssim \varepsilon^\alpha(1+t)^{-\frac{d+2k}{4}}\big(Y(t)^{1+\gamma_1}+Y(t)^2\big) \\ &\lesssim \varepsilon^\alpha(1+t)^{-\frac{d+2k}{4}}\big(1+Y(t)^2\big).\end{aligned} \tag{3.36}$$

Finally, we estimate $K_5$ using the arguments for $K_1$ and $K_3$. We therefore obtain

$$K_5 \leq (1+t)^{-\frac{d+2k}{4}}\left(C(\epsilon)+\epsilon Y(t)+\varepsilon^\alpha Y(t)^2\right). \tag{3.37}$$

Combining (3.35), (3.36) and (3.37), and adjusting $\epsilon$ once more, we conclude that

$$Q_4(t) \leq (1+t)^{-\frac{d+2k}{4}}\left(C(\epsilon)+\epsilon Y(t)+\varepsilon^\alpha Y(t)^2\right), \qquad 1\leq k\leq [d/2]+2.$$

Finally, we estimate the term $Q_6$.
Similarly, using the expression of $F_2$ and Minkowski's inequality, we obtain

$$\begin{aligned}Q_6(t) &\lesssim \left(\int_0^t\int_{\mathbb{R}^d} e^{-2c\lambda_\xi(t-\tau)}|\xi|^{2k+2}\left|\mathcal{F}\left(\frac{(\rho^\varepsilon-1)n^\varepsilon}{\rho^\varepsilon}\right)\right|^2 d\xi\, d\tau\right)^{\frac{1}{2}} \\ &\quad + \left(\int_0^t\int_{\mathbb{R}^d} e^{-2c\lambda_\xi(t-\tau)}|\xi|^{2k+2}\left|\mathcal{F}\left(\frac{z^\varepsilon m}{\rho\rho^\varepsilon}\right)\right|^2 d\xi\, d\tau\right)^{\frac{1}{2}}\end{aligned}$$

$$=: K_6 + K_7.$$

We first estimate $K_6$. A frequency decomposition gives

$$\begin{aligned}
K_6 \leq & \left( \int_0^t \int_{|\xi|\geq 1} e^{-2c\lambda_\xi(t-\tau)} |\xi|^{2k+2} \left| \mathcal{F}\left( \frac{(\rho^\varepsilon - 1)n^\varepsilon}{\rho^\varepsilon} \right) \right|^2 d\xi\, d\tau \right)^{\frac{1}{2}} \\
& + \left( \int_0^t \int_{|\xi|\leq 1} e^{-2c\lambda_\xi(t-\tau)} |\xi|^{2k+2} \left| \mathcal{F}\left( \frac{(\rho^\varepsilon - 1)n^\varepsilon}{\rho^\varepsilon} \right)(\tau) \right|^2 d\xi\, d\tau \right)^{\frac{1}{2}} \\
=: & K_6^1 + K_6^2.
\end{aligned}$$

For the high-frequency contribution, the Moser-type inequalities and Plancherel's identity imply

$$\begin{aligned}
K_6^1 \lesssim & \left( \int_0^t e^{-C(t-\tau)} \left\| D^{k+1} \left( \frac{(\rho^\varepsilon - 1)n^\varepsilon}{\rho^\varepsilon} \right) \right\|^2 d\tau \right)^{\frac{1}{2}} \\
\lesssim & \Bigg( \int_0^t e^{-C(t-\tau)} \Big( \left\| D^{k+1}(\rho, \varepsilon^\alpha z^\varepsilon) \right\|^2 \|n^\varepsilon\|_{L^\infty}^2 + \|(\rho - 1, \varepsilon^\alpha z^\varepsilon)\|_{L^\infty}^2 \|D^{k+1} n^\varepsilon\|^2 \\
& + \left\| D^{k+1}(\rho, \varepsilon^\alpha z^\varepsilon) \right\|^2 \|(\rho - 1, \varepsilon^\alpha z^\varepsilon)\|_{L^\infty}^2 \|n^\varepsilon\|_{L^\infty}^2 \Big) d\tau \Bigg)^{\frac{1}{2}} \\
\lesssim & \Bigg( \int_0^t e^{-C(t-\tau)} \Big( \|D^{k+1}\rho\|^2 \|n^\varepsilon\|_{L^\infty}^2 + \varepsilon^{2\alpha} \|D^k z^\varepsilon\|^{2\gamma_1} \|D^s z^\varepsilon\|^{2-2\gamma_1} \|n^\varepsilon\|_{L^\infty}^2 \\
& + \|(\rho - 1, \varepsilon^\alpha z^\varepsilon)\|_{L^\infty}^2 \|D^k n^\varepsilon\|^{2\gamma_1} \|D^s n^\varepsilon\|^{2-2\gamma_1} \Big) d\tau \Bigg)^{\frac{1}{2}}.
\end{aligned}$$

For $k = 1, \cdots, [\frac{d}{2}] + 1$, the assumption on $s$ in Theorem 1.2 gives $s \geq [\frac{d}{2}] + 4$, which implies

$$\gamma_1 \frac{d + 2k + 2}{4} + \frac{d}{2} \geq \frac{d + 2k}{4}.$$

It follows from (1.10) and assumption ($\mathrm{H}_2$) that

$$\begin{aligned}
K_6^1 \lesssim & \Bigg( \left( \int_0^{\frac{t}{2}} + \int_{\frac{t}{2}}^t \right) e^{-C(t-\tau)} \Big( (1+\tau)^{-\frac{d+2k+2}{2}} + Y(t)^{2\gamma_1} (1+\tau)^{-\gamma_1 \frac{d+2k+2}{2} - d} \\
& + \varepsilon^{2\alpha} Y(t)^{2+2\gamma_1} (1+\tau)^{-\gamma_1 \frac{d+2k+2}{2} - d - 1} \Big) d\tau \Bigg)^{\frac{1}{2}} \\
\lesssim & (1+t)^{-\frac{d+2k}{4}} (1 + Y(t)^{\gamma_1} + \varepsilon^\alpha Y(t)^{1+\gamma_1}).
\end{aligned} \tag{3.38}$$

When $k = [\frac{d}{2}] + 2$, we similarly obtain

$$\gamma_1 \frac{d + 2[\frac{d}{2}] + 4}{4} + \frac{d}{2} \geq \frac{d + 2[\frac{d}{2}] + 4}{4},$$

and

$$
\begin{aligned}
K_6^1 \lesssim & \left( \left( \int_0^{\frac{t}{2}} + \int_{\frac{t}{2}}^t \right) e^{-C(t-\tau)} \left( (1+\tau)^{-\frac{d+2[\frac{d}{2}]+6}{2}} + Y(t)^{2\gamma_1}(1+\tau)^{-\gamma_1 \frac{d+2[\frac{d}{2}]+4}{2}-d} \right.\right. \\
& \left.\left. + \varepsilon^{2\alpha} Y(t)^{2+2\gamma_1}(1+\tau)^{-\gamma_1 \frac{d+2[\frac{d}{2}]+4}{2}-d} \right) d\tau \right)^{\frac{1}{2}} \qquad (3.39) \\
\lesssim & (1+t)^{-\frac{d+2[\frac{d}{2}]+4}{4}} (1+Y(t)^{\gamma_1} + \varepsilon^\alpha Y(t)^{1+\gamma_1}).
\end{aligned}
$$

Fix $\gamma \in [\frac{3}{4}, 1)$. For the low-frequency part $K_6^2$ with $0 \le \tau \le \frac{t}{2}$, in view of (1.10) and assumption (H2), we obtain

$$
\begin{aligned}
& \left( \int_0^{\frac{t}{2}} \int_{|\xi| \le 1} e^{-\frac{C|\xi|^2(t-\tau)}{1+|\xi|^2}} |\xi|^{2k+2} \left| \mathcal{F}\left( \frac{(\rho^\varepsilon - 1)n^\varepsilon}{\rho^\varepsilon} \right) \right|^2 d\xi d\tau \right)^{\frac{1}{2}} \\
\lesssim & \left( \int_0^{\frac{t}{2}} \left\| \frac{(\rho^\varepsilon - 1)n^\varepsilon}{\rho^\varepsilon} \right\|_{L^1}^2 \int_{|\xi| \le 1} e^{-C|\xi|^2(1+t-\tau)} |\xi|^{2k} d\xi d\tau \right)^{\frac{1}{2}} \\
\lesssim & \left( \int_0^{\frac{t}{2}} \left( Y(t)^{2\gamma}(1+\tau)^{-\frac{d}{2}(1+\gamma)} + \varepsilon^{2\alpha} Y(t)^4 (1+\tau)^{-d-1} \right) (1+t-\tau)^{-\frac{d+2k}{2}} d\tau \right)^{\frac{1}{2}} \\
\lesssim & (1+t)^{-\frac{d+2k}{4}} (Y(t)^\gamma + \varepsilon^\alpha Y(t)^2).
\end{aligned}
$$

For the low-frequency part $K_6^2$ with $\frac{t}{2} \le \tau \le t$, we use the same $\gamma \in [\frac{3}{4}, 1)$. For all $k = 1, \cdots, [\frac{d}{2}] + 1$, we have

$$
d + \gamma \frac{d+2k+2}{2} \ge \frac{d+2k+2}{2}.
$$

For $k = [\frac{d}{2}] + 2$,

$$
d + \gamma \frac{d+2[\frac{d}{2}]+4}{2} \ge \frac{d+2[\frac{d}{2}]+6}{2}.
$$

Owing to (1.10), ($H_2$), the Moser-type inequalities and Young's inequality, we obtain

$$
\begin{aligned}
& \left( \int_{\frac{t}{2}}^t \int_{|\xi| \le 1} e^{-C|\xi|^2(t-\tau)} |\xi|^{2k+2} \left| \mathcal{F}\left( \frac{(\rho^\varepsilon - 1)n^\varepsilon}{\rho^\varepsilon} \right) \right|^2 d\xi \, d\tau \right)^{\frac{1}{2}} \\
\lesssim & \left( \int_{\frac{t}{2}}^t \int_{|\xi| \le 1} |\xi|^{2k} \left| \mathcal{F}\left( \frac{(\rho^\varepsilon - 1)n^\varepsilon}{\rho^\varepsilon} \right) \right|^2 d\xi \, d\tau \right)^{\frac{1}{2}} \\
\lesssim & \left( \int_{\frac{t}{2}}^t \left( \left\| D^k(\rho, \varepsilon^\alpha z^\varepsilon) \right\|^2 \|n^\varepsilon\|_{L^\infty}^2 + \|(\rho - 1, \varepsilon^\alpha z^\varepsilon)\|_{L^\infty}^2 \|D^k n^\varepsilon\|^2 \right.\right. \\
& \left.\left. + \left\| D^k(\rho, \varepsilon^\alpha z^\varepsilon) \right\|^2 \|\rho^\varepsilon - 1\|_{L^\infty}^2 \|n^\varepsilon\|_{L^\infty}^2 \right) d\tau \right)^{\frac{1}{2}} \\
\lesssim & (1+t)^{-\frac{d+2k}{4}} (Y(t)^\gamma + \varepsilon^\alpha Y(t)^2).
\end{aligned}
$$

Thus, we obtain the following estimate for $K_6^2$:

$$
K_6^2 \lesssim (1+t)^{-\frac{d+2k}{4}} (1 + Y(t)^\gamma + Y(t)^{\gamma_1} + \varepsilon^\alpha Y(t)^{1+\gamma_1} + \varepsilon^\alpha Y(t)^2). \qquad (3.40)
$$

Combining (3.38), (3.39), (3.40) and Young's inequality, we have

$$K_6 \lesssim (1+t)^{-\frac{d+2k}{4}}(C(\epsilon)+\epsilon Y(t)+\varepsilon^{\alpha}Y(t)^2).$$

The estimate for $K_7$ is obtained similarly. Hence, $(3.26)_3$ holds for $j=6$.

Substituting the estimates for $Q_0$-$Q_6$ in (3.26) into (3.25), we obtain

$$(1+t)^{\frac{d+2k}{4}}\|D^{k-1}(z^{\varepsilon}(t),n^{\varepsilon}(t))\|_1 \lesssim C(\epsilon)+\epsilon Y(t)+\varepsilon^{\alpha}Y(t)^2.$$

Summing the above inequality over $k=1,\dots,[\frac{d}{2}]+2$, noting that

$$\|(z^{\varepsilon},n^{\varepsilon})\|_{L^\infty} \lesssim \|(z^{\varepsilon},n^{\varepsilon})\|^{1-\theta}\|D^{[\frac{d}{2}]+1}(z^{\varepsilon},n^{\varepsilon})\|^{\theta}, \quad \theta=\frac{d}{2\,[d/2]+2},$$

and using (3.26), we obtain

$$\sum_{k=1}^{[\frac{d}{2}]+2}(1+t)^{\frac{d+2k}{4}}\|D^{k-1}(z^{\varepsilon}(t),n^{\varepsilon}(t))\|_1 \lesssim C(\epsilon)+\epsilon Y(t)+\varepsilon^{\alpha}Y(t)^2.$$

Since $\epsilon>0$ denotes a sufficiently small constant independent of $t,x$ and $\varepsilon$, the definition of $Y(t)$ gives (3.1).

The functional $Y(t)$ is continuous, and the initial-data assumptions bound $Y(0)$ uniformly in $\varepsilon$. Enlarging the constant term in (3.1) if necessary, we apply Strauss' lemma (Lemma A.1) to obtain $Y(t)\lesssim 1$ on $[0,T]$ for each fixed $\alpha\in(0,1]$ and sufficiently small $\varepsilon$. Since the constants are independent of $T$, the conclusion holds for all $t\geq 0$. □

# A Auxiliary lemmas

We collect two auxiliary lemmas used in the proof of Theorem 1.2. We first recall Strauss' lemma, which yields uniform bounds from nonlinear inequalities.

**Lemma A.1** ([29]). *Let $U(t)$ be a non-negative continuous function on $[0,T]$ satisfying the inequality*

$$U(t)\leq c_1+c_2U(t)^{\zeta}$$

*in some interval containing* 0*, where $c_1$ and $c_2$ are positive constants and $\zeta>1$. If $U(0)<c_1$ and*

$$c_1c_2^{\frac{1}{\zeta-1}}<\left(1-\frac{1}{\zeta}\right)\zeta^{-\frac{1}{\zeta-1}},$$

*then in the same interval*

$$U(t)<\frac{c_1}{1-\zeta^{-1}}.$$

We next establish a scalar comparison estimate that converts the Fourier energy inequality into an integral bound.

**Lemma A.2.** *Let $T > 0$, $a \geq 0$ and $C_0 > 0$. Let $A(t) \in L^1(0,T)$ and $B(t) \in L^2(0,T)$ be nonnegative. Suppose that $U(t) \in C([0,T])$ is nonnegative, $U(t)^2$ is absolutely continuous on $[0,T]$, and*

$$\frac{d}{dt}U(t)^2 + 2aU(t)^2 \leq C_0 A(t)U(t) + C_0 B(t)^2 \tag{A.1}$$

*for almost every $t \in (0,T)$. Define the nonnegative functions $Z(t)$ and $W(t)$ by*

$$Z(t) = e^{-at}U(0) + \frac{C_0}{2}\int_0^t e^{-a(t-\tau)}A(\tau)\,d\tau,$$
$$W(t)^2 = C_0\int_0^t e^{-2a(t-\tau)}B(\tau)^2\,d\tau.$$

*Then*

$$U(t) \leq Z(t) + W(t), \qquad 0 \leq t \leq T.$$

*Proof.* By definition,

$$\partial_t Z + aZ = \frac{C_0}{2}A, \qquad Z(0) = U(0),$$
$$\partial_t(W^2) + 2aW^2 = C_0 B^2, \qquad W(0) = 0.$$

For $\delta > 0$, define $W_\delta > 0$ by

$$W_\delta(t)^2 = W(t)^2 + \delta^2 e^{-2at}.$$

Then

$$\partial_t(W_\delta^2) + 2aW_\delta^2 = C_0 B^2, \qquad W_\delta(0) = \delta.$$

Consequently,

$$\partial_t(Z + W_\delta) + a(Z + W_\delta) = \frac{C_0}{2}A + \frac{C_0 B^2}{2W_\delta}.$$

On any interval where $U > Z + W_\delta$, we have $U > 0$. Dividing (A.1) by $2U$ gives

$$\partial_t U + aU \leq \frac{C_0}{2}A + \frac{C_0 B^2}{2U}.$$

Subtracting the preceding identity and using $U > Z + W_\delta \geq W_\delta > 0$, we obtain

$$\partial_t(U - Z - W_\delta) + a(U - Z - W_\delta) \leq \frac{C_0 B^2}{2}\left(\frac{1}{U} - \frac{1}{W_\delta}\right) \leq 0.$$

Hence, on such an interval,

$$\partial_t\big[e^{at}(U - Z - W_\delta)\big] \leq 0.$$

Since

$$U(0) - Z(0) - W_\delta(0) = -\delta < 0,$$

any nonempty connected component of

$$\{t \in [0,T] : U(t) > Z(t) + W_\delta(t)\}$$

would have a left endpoint $t_0 > 0$ at which

$$U(t_0) - Z(t_0) - W_\delta(t_0) = 0.$$

Integrating the preceding differential inequality from $t_0$ to a point in that component gives a contradiction. Therefore,

$$U(t) \leq Z(t) + W_\delta(t), \qquad 0 \leq t \leq T.$$

Letting $\delta \to 0^+$ proves the assertion. □

**Acknowledgment:** The research of Tao Wang is partially supported by the National Natural Science Foundation of China (Grants No. 12422109, 12426203, and 12371225). The research of Huijiang Zhao is supported by the National Natural Science Foundation of China (Grant No. 12221001, 12371225, and 12571242). Junhao Zhang is supported in part by the National Science Centre, Poland (Grant No. 2023/50/A/ST1/00447). Junhao Zhang also thanks Professor Agnieszka Świerczewska-Gwiazda for her support during his postdoctoral research and the University of Warsaw for its hospitality.

**Declaration of generative AI use:** During the preparation of this manuscript, the authors used generative AI tools to assist with English-language editing, checking mathematical arguments, and organizing relevant literature. The authors take full responsibility for the accuracy and integrity of the manuscript, including the mathematical results and references.

**Data availability statement:** No datasets were generated or analysed during the current study.

# References

[1] G. Boillat, Sur l'existence et la recherche d'équations de conservation supplémentaires pour les systèmes hyperboliques, C. R. Acad. Sci. Paris Sér. A **278** (1974), 909–912.

[2] T. Crin-Barat, S. Kawashima, J. Xu, The Cattaneo–Christov approximation of Fourier heat-conductive compressible fluids, Comm. Math. Phys. **406** (2025), no. 10, Paper No. 237, 70 pp.

[3] R. Danchin, Global existence in critical spaces for compressible Navier-Stokes equations, Invent. Math. **141** (2000), no. 3, 579–614.

[4] K. O. Friedrichs, P. D. Lax, Systems of conservation equations with a convex extension, Proc. Natl. Acad. Sci. U.S.A. **68** (1971), no. 8, 1686–1688.

[5] S. K. Godunov, An interesting class of quasi-linear systems (in Russian), Dokl. Akad. Nauk SSSR **139** (1961), no. 3, 521–523.

[6] Y. Guo, Y. Wang, Decay of dissipative equations and negative Sobolev spaces, Comm. Partial Differential Equations **37** (2012), no. 12, 2165–2208.

[7] D. Hoff, K. Zumbrun, Multi-dimensional diffusion waves for the Navier-Stokes equations of compressible flow, Indiana Univ. Math. J. **44** (1995), no. 2, 603–676.

[8] D. Hoff, K. Zumbrun, Pointwise decay estimates for multidimensional Navier-Stokes diffusion waves, Z. Angew. Math. Phys. **48** (1997), 597–614.

[9] Y. Hu, R. Racke, Compressible Navier-Stokes equations with revised Maxwell's law, J. Math. Fluid Mech. **19** (2017), no. 1, 77–90.

[10] Y. Hu, R. Racke, Global existence versus blow-up for multidimensional hyperbolized compressible Navier-Stokes equations, SIAM J. Math. Anal. **55** (2023), no. 5, 4788–4815.

[11] X. Huang, J. Li, Global well-posedness of classical solutions to the Cauchy problem of two-dimensional barotropic compressible Navier-Stokes system with vacuum and large initial data, SIAM J. Math. Anal. **54** (2022), no. 3, 3192–3214.

[12] Q. Jiu, Y. Wang, Z. Xin, Global classical solution to two-dimensional compressible Navier-Stokes equations with large data in $\mathbb{R}^2$, Phys. D **376/377** (2018), 180–194.

[13] Q. Jiu, Y. Wang, Z. Xin, Corrigendum to "Global classical solution to two-dimensional compressible Navier-Stokes equations with large data in $\mathbb{R}^2$" [Physica D 376/377 (2018) 180–194], Phys. D **414** (2020), Paper No. 132726.

[14] Y. Kagei, T. Kobayashi, Asymptotic behavior of solutions of the compressible Navier-Stokes equations on the half space, Arch. Ration. Mech. Anal. **177** (2005), no. 2, 231–330.

[15] T. Kato, The Cauchy problem for quasi-linear symmetric hyperbolic systems, Arch. Ration. Mech. Anal. **58** (1975), 181–205.

[16] P. D. Lax, *Hyperbolic Systems of Conservation Laws and the Mathematical Theory of Shock Waves*, CBMS-NSF Regional Conference Series in Applied Mathematics, vol. 11, SIAM, Philadelphia, 1973.

[17] F. Li, H. Tang, S. Zhang, Global well-posedness and large-time behavior of the compressible Navier-Stokes equations with hyperbolic heat conduction, J. Differential Equations **460** (2026), Paper No. 114111.

[18] T. P. Liu, W. Wang, The pointwise estimates of diffusion wave for the Navier-Stokes systems in odd multi-dimensions, Comm. Math. Phys. **196** (1998), no. 1, 145–173.

[19] A. Majda, *Compressible Fluid Flow and Systems of Conservation Laws in Several Space Variables*, Applied Mathematical Sciences, vol. 53, Springer-Verlag, New York, 1984.

[20] A. Matsumura, T. Nishida, The initial value problem for the equations of motion of compressible viscous and heat-conductive fluids, Proc. Japan Acad. Ser. A Math. Sci. **55** (1979), no. 9, 337–342.

[21] A. Matsumura, T. Nishida, The initial value problem for the equations of motion of viscous and heat-conductive gases, J. Math. Kyoto Univ. **20** (1980), no. 1, 67–104.

[22] J. Maxwell, On the dynamical theory of gases, Philos. Trans. R. Soc. Lond. **157** (1867), 49–88.

[23] J. Nash, Le problème de Cauchy pour les équations différentielles d'un fluide général, Bull. Soc. Math. France **90** (1962), 487–497.

[24] Y. J. Peng, Relaxed Euler systems and convergence to Navier-Stokes equations, Ann. Inst. H. Poincaré C Anal. Non Linéaire **38** (2021), no. 2, 369–401.

[25] Y. J. Peng, L. Zhao, Global convergence to compressible full Navier-Stokes equations by approximation with Oldroyd-type constitutive laws, J. Math. Fluid Mech. **24** (2022), no. 2, Paper No. 29, 17 pp.

[26] Y. J. Peng, L. Zhao, Global convergence rates from relaxed Euler equations to Navier-Stokes equations with Oldroyd-type constitutive laws, Nonlinearity **37** (2024), no. 9, Paper No. 095032, 26 pp.

[27] J. Serrin, On the uniqueness of compressible fluid motions, Arch. Ration. Mech. Anal. **3** (1959), 271–288.

[28] T. C. Sideris, Formation of singularities in solutions to nonlinear hyperbolic equations, Arch. Ration. Mech. Anal. **86** (1984), no. 4, 369–381.

[29] W. Strauss, Decay and asymptotics for $\Box u = F(u)$, J. Funct. Anal. **2** (1968), 409–457.

[30] H. Tang, S. Zhang, W. Zou, Decay of the compressible Navier-Stokes equations with hyperbolic heat conduction, J. Differential Equations **388** (2024), 1–33.

[31] V. A. Vaigant, A. V. Kazhikhov, On existence of global solutions to the two-dimensional Navier-Stokes equations for a compressible viscous fluid, Sib. Math. J. **36** (1995), no. 6, 1108–1141.

[32] N. Wang, Y. Hu, Blowup of solutions for compressible Navier-Stokes equations with revised Maxwell's law, Appl. Math. Lett. **103** (2020), Paper No. 106221.

[33] W. A. Yong, Newtonian limit of Maxwell fluid flows, Arch. Ration. Mech. Anal. **214** (2014), no. 3, 913–922.